\documentclass[10pt,a4paper,dk,reqno]{amsart}
\let\amsmarkboth\markboth    
\usepackage{amsmath,amsfonts,amssymb,amsthm,amscd}
\usepackage[english]{babel}
\usepackage[mathscr]{eucal}
\usepackage{graphicx}
\usepackage{appendix}
\usepackage{color}

\makeatletter
\let\markboth\amsmarkboth
\bbl@redefine#1#2{%
   \def\bbl@arg{#1}%
   \ifx\bbl@arg\@empty
     \toks@{}%
   \else
     \toks@{\noexpand\foreignlanguage{%
              \languagename}{%
              \noexpand\bbl@restore@actives#1}}%
   \fi
   \def\bbl@arg{#2}%
   \ifx\bbl@arg\@empty
     \toks8{}%
   \else
     \toks8{\noexpand\foreignlanguage{%
              \languagename}{%
              \noexpand\bbl@restore@actives#2}}%
   \fi
   \edef\bbl@tempa{\the\toks@}%
   \edef\bbl@tempb{\the\toks8}%
   \protected@edef\bbl@tempa{%
     \noexpand\org@markboth{\bbl@tempa}{\bbl@tempb}}%
   \bbl@tempa
}
\DeclareRobustCommand*\ams@disablelinebreak{\def\\{ \ignorespaces}}
\def\maketitle{\par
   \@topnum\z@ %
   \@setcopyright
   \thispagestyle{firstpage}%
   \uppercasenonmath\shorttitle
   \ifx\@empty\shortauthors \let\shortauthors\shorttitle
   \else \andify\shortauthors
   \fi
   \@maketitle@hook
   \begingroup
   \@maketitle
   \toks@\@xp{\shortauthors}\@temptokena\@xp{\shorttitle}%
   \protected@edef\@tempa{%
     \@nx\markboth{\ams@disablelinebreak
       \@nx\MakeUppercase{\the\toks@}}{\the\@temptokena}}%
   \@tempa
   \endgroup
   \c@footnote\z@
   \@cleartopmattertags
}

\makeatother

\numberwithin{equation}{section}

\newtheorem{remark}[]{Remark}
\theoremstyle{definition}

\renewcommand{\P}{\Psi}

\newcommand{\C}{\mathbb{C}}
\newcommand{\R}{\mathbb{R}}

\newcommand{\eps}{\varepsilon}
\newcommand{\dsp}{\displaystyle}

\title[]{Evolution, interaction and collisions\\ of vortex filaments}

\author[]{Valeria Banica}
\address[V. Banica]{Laboratoire Analyse et probabilit\'es (EA                  
2172), D\'epartement de Math\'ematiques, Universit\'e d'Evry, 23 Bd. de       
France, 91037 Evry, France, Valeria.Banica@univ-evry.fr}

\author[]{Evelyne Miot}
\address[E. Miot]{Laboratoire de Math\'ematiques, Universit\'e Paris-Sud 11, B\^at. 425, 91405 Orsay, France, Evelyne.Miot@math.u-psud.fr}

\begin{document}
\maketitle

\begin{abstract}
Several progresses have been done very recently on models for the dynamics of one or more vortex filaments in 3-D fluids. In this article we survey the recent and previous  results in this topic. We finally present a collection of new situations of filaments collapse. 
\end{abstract}
 \tableofcontents

\section{Introduction}

The purpose of this work is to describe some aspects of the dynamics of homogeneous three-dimensional incompressible fluids.  The evolution of such fluids is
governed by the Euler equations
\begin{equation}
\label{eq:Euler}
\partial_t \mathbf{\omega}+v\cdot \nabla \omega=\omega\cdot \nabla v,\quad \text{div}\,v=0,
\end{equation}
where $v:\R\times \R^3\to \R^3$  and $\omega=\text{curl}(v):\R\times \R^3\to \R^3$ denote the velocity and the vorticity of the fluid.

We shall focus on fluids in which remarkable structures, called vortex filaments, are present. A vortex filament is a  very thin tube in which the vorticity is sharply concentrated.  It can be asymptotically assimilated to a curve $\chi(t,s)$ in $\R^3$, where $t\in \R$ denotes the time and $s\in \R$ the arc-length parameter. According to a suitable formal derivation, which will be described in the next section, the asymptotic motion law of the curve is the binormal flow 
\begin{equation}
\tag{BF}\label{eq:BF}
\chi_t=\partial_s \chi \wedge\partial_{s}^2\chi.
\end{equation}
Equation \eqref{eq:BF} exhibits a rich variety of motions, which will be presented in Section \ref{Sec:one-filament} below. Particular attention will be payed to motions generating a singularity in finite time.

\medskip

On the other hand,  special solutions of the binormal flow are given by the infinite straight filaments. Up to a rotation, we can assume that there are parallel to the $e_3$--axis so that they can be parametrized as
$\chi(t,\sigma)=(X_0,\sigma)$ for some fixed $X_0\in \R^2$, where $\sigma\in \R$ can be also seen as an arc-length parameter. Of interest is the evolution of perturbations of such straight filaments, namely the nearly parallel vortex filaments, which are parametrized by
$\chi(t,\sigma)=(X_0+\Psi(t,\sigma),\sigma)$ where $\Psi(t,\sigma)\in \R^2\simeq \C$ is small. It turns out that linearization about \eqref{eq:BF} leads to the linear one-dimensional Schr\"odinger equation for $\Psi$
\begin{equation}
\label{eq:filament}
i\partial_t \Psi+\partial_\sigma^2\Psi=0.
\end{equation}
Next,  we wish to investigate the dynamics of several vortex filaments that are all nearly parallel to the same $e_3$--axis. We consider a collection of $N\geq 2$ filaments, which are represented by the complex-valued functions $\Psi_j(t,\sigma)\in \C$, with $t\in \R$ and $\sigma\in \R$. 
In 1995, Klein, Majda and Damodaran \cite{KlMaDa} derived, in a suitable asymptotic regime, a system of coupled equations for the filament positions $\Psi_j$
 \begin{equation}
 \label{syst:filaments}
\displaystyle 
i\partial_t \Psi_j+\alpha_j\Gamma_j \partial_\sigma^2\Psi_j+\sum_{k\neq j} \Gamma_k \frac{\Psi_j-\Psi_k}{|\Psi_j-\Psi_k|^2}=0,\quad 1\leq j\leq N. 
\end{equation}
Here $\Gamma_j\in \R$ corresponds to the circulation of the $j$--th filament, and the parameter $\alpha_j\in \R$ is related to the core structure of the filament. In particular, in the case where the vortex filaments are all exactly parallel, i.e. if $\Psi_j(t,\sigma)=X_j(t)$ does not depend on $\sigma$,  the previous system reduces to a two-dimensional system of ordinary differential equations for the positions $X_j(t)$
\begin{equation}
 \label{syst:point-vortex}
i\frac{d X_j}{dt}+\sum_{k\neq j} \Gamma_k \frac{X_j-X_k}{|X_j-X_k|^2}=0,\quad 1\leq j\leq N. 
\end{equation}
This system is called point vortex system or Kirchhoff law, and it has been intensively studied in the literature. 

We observe that the system \eqref{syst:filaments} combines two different aspects of the dynamics. On the one hand, the linearized self-induced motion of each filament --already expressed in the  equation \eqref{eq:filament}-- is represented through the Schr\"odinger operator. On the other hand, the interaction of the $j$--th filament with the other filaments is represented by the potential field $\sum_{k\neq j} \Gamma_k (\Psi_j-\Psi_k)/|\Psi_j-\Psi_k|^2$, which corresponds to the velocity generated by the other filaments\footnote{According to the Biot-Savart law, see Section \ref{Sec:one-filament}.}. For straight filaments there is no self-induced motion, which means that each filament moves only with the velocity induced by the other filaments according to \eqref{syst:point-vortex} -- exactly as for planar vortices, see Section \ref{Sec:point} hereafter.  The fact that both effects are taken into account in \eqref{syst:filaments} is really due to the  asymptotic conditions  chosen for the derivation, which relate the wavelength of perturbations,  the core sizes and the separation distances between the vortex tubes in a particular way.

Since the system \eqref{syst:filaments} is thought to be a simplified model for the dynamics of filaments, a first natural issue is the existence and uniqueness of solutions. The potential fields are not well-defined at points when two or more filaments collide in finite time.  One possibility  is to consider weak distributional solutions, allowing possibly for exceptional collisions, and for which the potential is well-defined almost-everywhere (see Section \ref{sec:filaments}). But then uniqueness is not achieved in this class. Therefore, we will focus on a class of  stronger solutions, for which the filaments are well-separated up to  the largest time of existence of the system. Hence the largest time of existence corresponds to the first collision time in this framework. It should be mentioned that,  even for the simpler situation of straight filaments (the point vortex system), collisions in finite time are known to occur for some initial configurations (see Section \ref{Sec:point}). Therefore we cannot hope for global existence in general. We shall determine sufficient conditions on the initial collection of vortex filaments that lead to  large time or global existence, and we shall also study the case of finite-time collapses.

\medskip

The remainder of this work is organized as follows. In Section \ref{Sec:one-filament}  we present the formal derivation of the binormal flow \eqref{eq:BF} for the 
motion of one single filament and we 
describe some of its remarkable properties. In particular we review the recent results on singularities for self-similar and almost self-similar solutions of  \eqref{eq:BF}, as well as  for their perturbations. Section \ref{Sec:point} contains a brief state-of-the-art concerning the point vortex system \eqref{syst:point-vortex}, as an introduction to the more complicated dynamics of almost parallel vortex filaments. Then, Section \ref{sec:filaments} gives an overview of existence and uniqueness results concerning the system \eqref{syst:filaments}. In particular, we present some recent results from \cite{BaMi} on symmetric configurations of vortex filaments. We conclude by providing several and simple examples of finite-time collapses of
vortex filaments. The collapses are of different nature, and some of them are new.

\medskip

\noindent \textbf{Notations}.
In the following we will systematically identify complex numbers and real vectors.  In particular, introducing the matrix $\mathbf{J}=\begin{pmatrix} 0&-1 \\1 & 0\end{pmatrix}$
we shall identify the $\R^2$-vector $\mathbf{J}X$ and the complex $iX$.

\medskip

\noindent \textbf{Acknowledgements}. The first author is partially supported by the ANR project ``R.A.S.''. Both authors are grateful to Vincent Torri for his support in using Scilab.

\section{One filament}\label{Sec:one-filament}

In this section we briefly present classical facts about the local 
induction approximation and about the binormal flow, as well as recent
  results of singularity formation in finite time. For parts of this
description the reader may also consult the book \cite{MaBe} and the 
surveys \cite{ponce} and \cite{BaVe4}.

\subsection{The Local Induction Approximation}
We shall start with a short description of the way the binormal flow equation was find as a model for vortex filaments dynamics in a three-dimensional ideal incompressible fluid. 
A fluid with a vortex filament is a fluid where the vorticity $\omega(t)$ is a singular measure supported along a curve $\chi(t)$, which is parametrized by arc-length parameter $s$, and with density  $\Gamma \partial_s \chi(t)$. Here $\Gamma$ denotes the constant circulation along the filament. As a consequence of Kelvin's circulation law vortex filaments move with the flow, so in order to get the dynamics of the vortex curve $\chi(t)$ we have to compute the velocity of the fluid near the curve. The three-dimensional Biot-Savart kernel is 
$\frac{1}{4\pi|x|}$ so the velocity of the fluid is
$$v(t,x)=-\nabla\wedge\left(\frac{1}{4\pi|\cdot|}*\omega(t,\cdot)\right)(x)=-\frac{\Gamma}{4\pi}\int_{-\infty}^\infty\frac{x-\chi(t,s)}{|x-\chi(t,s)|^3}\wedge\partial_s\chi(t,s)ds.$$
We fix $t$, we suppose without loss of generality that $\chi(t,0)=(0,0,0),\partial_s\chi(t,0)=(0,0,1)$ and we get focused on what is happening for the velocity of the fluid near $(0,0,0)$. For this purpose, two localizations are done, one is that $\chi(t,s)$ will be approximated by a Taylor development of order 2 near $s=0$ and the other is that the integral has to be considered only locally around $s=0$, on a segment $[-L,L]$. We obtain that $v(t,(\varepsilon,0,0))$ is approximated by 
$$-\frac{\Gamma}{4\pi}\frac{(\varepsilon,0,0)\wedge (0,0,1)}{\varepsilon^2}\int_{-\frac L\varepsilon}^{\frac L\varepsilon}\frac{ds}{|1+s^2|^\frac 32} + \frac{\Gamma}{4\pi}\frac{\partial_s \chi(t,0)\wedge \partial_s^2 \chi(t,0)}{2}\int_{-\frac L\varepsilon}^{\frac L\varepsilon}\frac{s^2ds}{|1+s^2|^\frac 32}.$$
When $\eps\to 0$, the first term is diverging but in the same way the velocity does around one straight vortex filament. One straight vortex filament remains still in a fluid, so this first contribution should not influence the dynamics of one single vortex filament $\chi(t)$. Finally, even after localization, the integral in the last term is still diverging as $|\ln \varepsilon|$. A time rescaling allow to resorb this growth and then the binormal flow
\begin{equation}\label{bf}
\partial_t\chi=\partial_s\chi\wedge\partial_s^2\chi
\end{equation}
appears as an approximation for the dynamics of vortex filaments, called Local Induction Approximation. It was discovered by Da Rios \cite{DaR} in 1906 and re-lightened by Arms and Hama \cite{ArHa} in 1965. However the rigorous derivation of \eqref{bf} from the Euler equations \eqref{eq:Euler}, and in fact even the question of the stability of the vortex filaments, are open.

Being a crude model, the binormal flow has advantages and disadvantages. The positive points are the fact that special solutions as straight lines, translating circles, helices are on the one hand solutions of the binormal flow and on the other hand these dynamics are observed in fluid dynamics. This way other special solutions of the binormal flow may give some intuition for particular type of vortex filaments in fluids. A successful example in this direction are the traveling type solutions of the binormal flow found in \cite{Ha} and displayed in an fluid experiment \cite{HoBr}. Also, the binormal flow is accessible for numerics. One of the main issues of the L.I.A. is that closed curves conserve their lengths under the binormal flow which is not the case for closed vortex filaments in fluids. 

\medskip

The binormal curvature flow also appears in the context of superfluids that are governed by the Gross-Pitaevskii equation
\begin{equation}\label{eq:GP-N}\tag{GP}
i\partial_t \Psi+\Delta \Psi+\Psi(1-|\Psi|^2)=0.\end{equation}
More precisely, \eqref{eq:BF} is conjectured to govern the asymptotic dynamics of  codimension-$2$ submanifolds around which the Jacobians of the solutions to (a rescaled version of) \eqref{eq:GP-N}  concentrate. This was proved in any dimension $N$ in the special case of the $N-2$-dimensional sphere \cite{Je}.  In a parallel setting, mean curvature flow governs the vortex dynamics for parabolic Ginzburg-Landau equation \cite{BeOrSm}.  The fact that vortex filaments exhibit a common asymptotic motion law in fluids and superfluids results from a fundamental analogy between \eqref{eq:GP-N} and the Euler equations in such regimes.  
More precisely, the Madelung \cite{Ma} transform $\Psi=\sqrt{\rho}\exp(i\varphi)$, as long as $\Psi$ does not vanish, yields a compressible Euler type system for the variables $\rho$ and $v=\nabla \varphi$. It is called hydrodynamical form of \eqref{eq:GP-N} and it reduces formally to the incompressible Euler equations for $v$ in the above-mentioned singular limit of \eqref{eq:GP-N}.   Further details may be found in the recent survey \cite{CaDaSa} and in references quoted therein. 

\subsection{The binormal flow} 
	The binormal flow is a rich geometric equation. It is a completely integrable system.  In the following we shall denote curvature and torsion of a solution of \eqref{bf} by $c(t,s)$ and $\tau(t,s)$.  Some remarkable  conserved quantities of the binormal flow are the kinetic energy $\int c^2(t,s)ds$, helicity $\int c^2\tau(t,s)ds$, linear momentum $\int\chi\wedge\partial_s\chi(t,s)ds$, angular momentum $\int\chi\wedge(\chi\wedge\partial_s\chi)(t,s)ds$.  It is a reversible in time equation, invariant with respect to translations and rotations. The tangent vector $T(t,s)$ of a solution of the binormal flow satisfies the Schr\"odinger map equation
\begin{equation}\label{schmap}\partial_tT=T\wedge \partial_s^2T\end{equation}
	which plays an important role in ferromagnetism as a simplification of the Landau-Lifschitz equation.\par
	Let us present now the link with the Schr\"odinger equation. As mentioned in the introduction, the linear Schr\"odinger equation can be obtained easily as a rough approximation of the binormal flow. This can be seen  by plugging in \eqref{bf} the ansatz 
	$$\chi(t,s)=(\varepsilon x(t,s),\varepsilon y(t,s),s).$$ To leading order when $\varepsilon\to 0$, we obtain that $x+iy$ is a solution of the linear Schr\"odinger equation. It was in 1972 that Hasimoto \cite{Ha} has discovered the link between the binormal flow and the cubic one-dimensional nonlinear Schr\"odinger equation. The so-called filament function constructed from curvature and torsion of a solution of \eqref{bf},
$$\psi(t,s)=c(t,s)e^{i\int_0^s\tau(t,r)dr},$$
 satisfies the 1-D cubic NLS
$$i\partial_t\psi+\partial_s^2\psi+\frac 12\left(|\psi|^2-A(t)\right)\psi=0,$$
with $A(t)$ in terms of curvature and torsion $(c,\tau)(t,0)$. More precisely,  
$$A(t)=\left(2\frac{\partial_s c^2-c\tau^2}{c}+c^2\right)(t,0),$$ which might be seen as a compatibility condition. Let us notice that the Hasimoto transform can be seen as an inverse of Madelung's transform - actually curvature and torsion satisfy an Euler with quantum pressure type of system.  One can remark that the non-vanishing of the curvature appears as a constraint. Nevertheless, Koiso gave in \cite{Ko} a way to avoid this issue by using, instead of Frenet frames, parallel frames $(T,e_1,e_2)$ that satisfy
$$
\partial_s\left(\begin{array}{c}
T\\e_1\\e_2
\end{array}\right)=
\left(\begin{array}{ccc}
0 & \alpha & \beta \\ -\alpha & 0 & 0 \\ -\beta &  0 & 0 
\end{array}\right)
\left(\begin{array}{c}
T\\e_1\\e_2
\end{array}\right)\,\,\,.$$
Then $\alpha+i\beta$ is a solution of 1-D cubic NLS. For details of this kind of Hasimoto transform we refer to Appendix A of \cite{ponce}. As a conclusion, up to a change of phase, the filament function satisfies the 1-D cubic NLS. 

By using the Hasimoto transform the first local well-posedness results were obtained for the binormal flow - curvature and torsion were considered in high order Sobolev spaces, 
see \cite{Ha, Ko,FuMi}. \par
A recent result of local well-posedness for the binormal flow, for less regular closed curves, was obtained by Jerrard and Smets \cite{JeSm1,JeSm2} by considering a weak version of the binormal flow. Their method acts at the level of the tangent of the curve and it does not use the Hasimoto transform. Also, they have derived this way new results for the Schr\"odinger map equation \eqref{schmap}.

As far as we have seen, there exists a number of intimate relations between the evolution of filaments in fluids and dispersive equations. We end this subsection by discussing other similar connections.
In order to improve the validity of the binormal flow as a model for vortex filaments dynamics, an extended version of \eqref{bf} was considered by  Fukumoto and Miyazaki \cite{FuMi}, which includes an axial flow along the filament,
\begin{equation*}
\partial_t\chi=\partial_s\chi\wedge\partial_s^2\chi+\alpha\left(\partial_s^3\chi+\frac 32\partial_s^2\chi\wedge(\partial_s\chi\wedge\partial_s^2\chi)\right).
\end{equation*}
Here $\alpha$ denotes the magnitude of the axial flow. The Hasimoto transform leads then to the so-called Hirota equation
$$i\partial_t\psi+\partial_s^2\psi+\frac 12|\psi|^2\psi-i\alpha\,\partial_s^3\psi+\frac 32|\psi|^2\partial_s\psi=0,$$
which is in link both with NLS and the modified KdV equation (see \cite{La} for the initial derivation of this equation). It actually turned out that the modified KdV equation appears as a model for the evolution of the curvature of a vortex patch in a two-dimensional fluid \cite{GoPe1, GoPe2, ZaHuRo}, see also  \cite{PeVe} for more details and for self-similar solutions of this model.

\subsection{Singularity formation scenarios for the binormal flow}

	An interesting class of solutions that is not covered by the previous well-posedness results are the self-similar solutions of the binormal flow and their small and smooth perturbations. Since solutions of \eqref{bf} are considered to be parametrized by arc-length, self-similar solutions are of type
	$$\chi(t,s)=\sqrt {t} G\left(\frac{s}{\sqrt{t}}\right).$$
It is then easy to compute that they form a family of solutions $\chi_a$ indexed by a real parameter $a>0$ and $(c_a,\tau_a)(t,s)=\left(\frac{a}{\sqrt{t}},\frac{s}{2t}\right)$ (\cite{La}). They have been studied in the last decades as a frictionless model for the evolution after a line-line reconnection in the absence of counterflow for superfluid 4He and in ferromagnetism \cite{Sc, Bu,Li1, LaRuTh, LaDa}. They also appear, this time for finite length filaments with appropriate boundary conditions, in some models for the fiber architecture of aortic valve leaflets, see \cite{PeMc}. Their rigorous mathematical study was achieved by Gutierrez, Rivas and Vega, who gave in \cite{GuRiVe} a description of the profile and in particular established that 
	$$\chi_a(0,s)=sA^+  \mathbb{I}_{[0,\infty)}(s)-sA^-  \mathbb{I}_{(-\infty,0]}(s),$$
with $A^\pm\in\mathbb S^2$ distinct, non-opposite and $\sin\frac{(\widehat{A^+,-A^-})}{2}=e^{-\frac{a^2}{2}}.$ So self-similar solutions of the binormal flow are smooth infinite curves that generate a corner in finite time. In fluid mechanics, their dynamics is close to the one of the delta-wing vortex. In \cite{Patxi} numerical simulations for self-similar solutions of the binormal flow were shown to be in correlation with the physical experiment.

 Another family of \eqref{bf} generating a singularity in finite time at the level of the Frenet frame are the almost self-similar solutions
 $$\chi(t,s)=e^{\mathcal A\log t}\sqrt {t} G\left(\frac{s}{\sqrt{t}}\right),$$
 where $\mathcal A$ is a real $3\times3$ antisymmetric matrix, see  \cite{Li1,Li2}. Guti\'errez and Vega \cite{GuVe1} proved that
 $$\chi_a(0,s)=se^{\mathcal A\log |s|}A^+  \mathbb{I}_{[0,\infty)}(s)-se^{\mathcal A\log |s|}A^-  \mathbb{I}_{(-\infty,0]}(s).$$
 The curves at times $t=0$ are spirals, some of them with a singularity point, see \cite{web} for some numerical simulations. Both similar and almost self-similar solutions have self-similar curvature and torsion, so the filament function is also of type $\psi(t,s)=\frac {1}{\sqrt {t}}G\left(\frac{s}{\sqrt {t}}\right)$ and has to solve
 \begin{equation}\label{Hasi}
 i\partial_t\psi+\partial_s^2\psi+\frac 12\left(|\psi|^2-\frac {a}{t}\right)\psi=0.
 \end{equation}
Therefore the data has to be homogeneous of degree $-1$. The case of Dirac distribution leads to the self-similar solutions of \eqref{bf}, and the principal value distribution to the almost self-similar solutions respectively.

The study of small and smooth perturbations of self-similar solutions was initiated by Banica and Vega \cite{BaVe0, BaVe1, BaVe2, BaVe3} from the point of view of the stability of the formation of the singularity. The starting point of the analysis is the Hasimoto transform: the filament function of the self-similar solutions $\chi_a$ is $\psi_a(t,x)=\frac{a}{\sqrt{t}}\,e^{i\frac{x^2}{4t}},$
solution of \eqref{Hasi}. By the pseudo-conformal transformation 
$$\psi(t,x)=\mathcal {T}(u+a)(t,x)=\frac{1}{\sqrt{t}}\,e^{i\frac{x^2}{4t}}\,\overline{u+a}\left(\frac 1t,\frac xt\right),$$
understanding the behavior at time $0$ for perturbations of $\psi_a$, solutions of \eqref{Hasi}, is equivalent to understanding the large time behavior for the equation
\begin{equation}\label{Hasi-2}
i\partial_tu+\partial_s^2u+\frac 1{2t}\left(|u+a|^2-a^2\right)(u+a)=0.
 \end{equation} A scattering theory has been developed in \cite{BaVe1, BaVe2} for this problem, by strongly relying on Fourier analysis. The informations on \eqref{Hasi-2} have been transported at the level of the binormal flow by treating differently the regions $ts^2\lesssim 1$ and $1\lesssim ts^2$. 
 In order to discard the linear non-oscillating term, the asymptotic profile for $u(t)$ has to be modified into the long-range one $$e^{i \frac{a^2}{2}\log t}\,e^{it\partial_x^2}\,u_+(x).$$
 
In  \cite{BaVe1} it has been shown that wave operators exist, with a strong $H^s$ decay, for asymptotic profiles $u_+\in \dot{H}^{-2}\cap H^s\cap W^{s,1}$. By adding the extra-conditions $x^2u_+\in L^2$ and $s\geq 3$ the existence of solutions of \eqref{bf} that develop a corner-type singularity at $t=0$ has been derived. This insures that the formation of singularity in finite time appearing in the self-similar case is not an isolated phenomenon.  A similar result for the almost self-similar solutions of \eqref{bf} has been proved in \cite{GuVe2}.

Asymptotic completeness with loss of regularity was then proved in  \cite{BaVe2}, in spaces defined in terms of the behavior of the low Fourier modes. At the level of the binormal flow, this information has been translated into the fact that all small and regular perturbations at time $t_0\neq 0$ of a self-similar solution $\chi_a$ still yield, after evolving by the binormal flow, a singularity at time $t=0$. 

Finally, in \cite{BaVe3} we give a geometrical description of the curves at the singularity time by getting a control of the evolution of weighted norms for equation  \eqref{Hasi-2}. We prove the stability of the selfsimilar dynamics of small pertubations of a given selfsimilar solution, and in particular the fact that the angle of $\chi_a$ is recovered at time $t=0$.

\section{The point vortex system}\label{section:vortexpoints}\label{Sec:point}

This section is devoted to a short presentation of the point vortex system \eqref{syst:point-vortex}. The reader may find the subsequent results as well as many additional details in the surveys \cite{Aref83, Neil87, ANSTV} and in references quoted therein.

The point vortex system arises as an asymptotic motion law for point singularities in two-dimensional incompressible fluids, which are governed by the incompressible Euler equations, or superfluids, governed by the Gross-Pitaevskii equation \eqref{eq:GP-N}.  In contrast with the three-dimensional case, the persistence of point vortices and the validity of the point vortex model have been established rigorously for well-prepared data,  at least up to the first collision time between the vortices (see \cite{Tu, MaPu-loc} for the Euler equations and, e.g., \cite{CoJe, BeJeSm} for the Gross-Pitaevskii equation).

\medskip

System \eqref{syst:point-vortex} possesses an Hamiltonian structure, since it has the form
\begin{equation*}
\displaystyle \Gamma_j\frac{d}{dt} X_j(t)+\mathbf{J}\nabla_{X_j} H(X_1,\ldots,X_N)=0,\quad  1\leq j\leq N,
\end{equation*}
with the Hamiltonian $H$ defined by
\begin{equation*}
H(X_1,\ldots,X_N)=-\frac{1}{2}\sum_{j\neq k} \Gamma_j \Gamma_k \ln |X_j-X_k|.
\end{equation*}
As a consequence, we have  the conservation of $H$
$$H(X_1(t),\ldots,X_N(t))=H(X_1(0),\ldots,X_N(0)).$$ 	Moreover,  the point
$$C=\sum_{j=1}^N \Gamma_j X_j(t)=\sum_{j=1}^N \Gamma_j X_j(0)$$
and the angular momentum
$$ \sum_{j=1}^N \Gamma_j |X_j(t)|^2=\sum_{j=1}^N \Gamma_j |X_j(0)|^2$$
are preserved in time. If $\sum_j\Gamma_j\neq 0$, we define the constant center of vorticity
$$c=\:(\sum_{j=1}^N \Gamma_j)^{-1}C.$$
It follows that the quantity
$$ T=\sum_{j, k} \Gamma_j \Gamma_k|X_j(t)-X_k(t)|^2=\sum_{j, k} \Gamma_j 
\Gamma_k|X_j(0)-X_k(0)|^2$$ is constant as well. A direct consequence of the previous conservation laws is that,  in case of circulations $\Gamma_j$ all having the same sign, no collision between the vortices can occur in finite time and the $N$-vortex system admits a unique and global solution.

\medskip

It is well-known that the point vortex system admits particular solutions that keep the geometric shape of the polygon formed by the vortex positions.  To describe these solutions  we shall adopt  the following definitions from, e.g.,  \cite{Neil87}: we say that the configuration $(X_j(t))_j$ is 

\textbullet$\:$  a fixed equilibrium if $X_j(t)=X_j(0)$, $\forall j=1,\ldots,N$;

\textbullet$\:$ rigidly translating if, for all $t$, there exists $V\in \C^\ast$ such that $\dot{X}_j(t)=V$, $\forall j=1,\ldots,N$;

\textbullet$\:$ a relative equilibrium if there exists some $X_0\in \C$ and for all $t$, there exists $\omega \in \R^\ast$ such that
$ \dot{X}_j(t)=i\omega(t) (X_j(t)- X_0)$, $\forall j=1,\ldots,N$;

\textbullet$\:$ a {self-similar} (or homographic) motion if there exists some $X_0\in \C$ and for all $t$, there exists $\omega \in \C$, with $\Im m(\omega)\neq 0$, such that
$ \dot{X}_j(t)=i\omega (X_j(t)- X_0)$,  $\forall j=1,\ldots,N$.

Note that in the first three cases both the shape and the size of the polygon formed by the points  are preserved, since the motion is a translation in the second case and a rotation--
with angular velocity $ \omega(t)$-- in the second case. In the case of homographic motion the configuration at time $t$  is obtained from the initial
one by a rotation and a dilation. More precisely,  it must satisfy $X_j(t)-X_0=(X_j(0)-X_0)\sqrt{1-t/t_0}\exp(i\lambda \ln|1-t/t_0|)$, where the reals $t_0$ and $\lambda$ are determined in terms of $\omega$. When $t_0>0$, which actually corresponds to the case $\Im m(\omega(t))> 0$,  the configuration shrinks  and collapses into the center of vorticity $X_0$ at time $t_0$. Otherwise we have $t_0<0$, i.e. $\Im m(\omega(t))< 0$ and the configuration expands in time. 

 \medskip

For $N=2$ the distance $|X_1(t)-X_2(t)|$ remains constant therefore the solution is global. The motion is a uniform rotation about the center of vorticity if $\Gamma_1+\Gamma_2\neq 0$, and a rigid uniform translation otherwise.

\medskip

For $N=3$ the motion is integrable and has been widely investigated, see e.g. \cite{Grobli, Aref79, ANSTV, Aref10}. An important part of the analysis was initiated in 1877 by Gr\"obli \cite{Grobli} in his PhD thesis. In particular Gr\"obli noticed that the $3$-vortex system  can be entirely formulated in terms of the three distances between the vortices. The motion can then be analyzed qualitatively through the use of trilinear coordinates and phase diagrams, see the work of Aref \cite{Aref79}. There exists a large variety of regimes: fixed equilibria (if  $\Gamma_1\Gamma_2+\Gamma_1\Gamma_3+\Gamma_2\Gamma_3=0$), rigidly translating equilateral triangles (if $\Gamma_1+\Gamma_2+\Gamma_3= 0$ and $T\neq 0$), uniformly rotating equilateral triangles (if $\Gamma_1+\Gamma_2+\Gamma_3\neq 0$ and $T\neq 0$), collinear rotating or fixed configurations. Also,  self-similar (shrinking or expanding) collapsing configurations exist \cite{Grobli, Aref79, MP, Aref10}. It turns out that the conditions
 $T=0$ and $\Gamma_1\Gamma_2+\Gamma_1\Gamma_3+\Gamma_2\Gamma_3=0$ are both necessary and sufficient to obtain a self-similar expanding or shrinking motion\footnote{If  the initial configuration is not a fixed equilibrium.}, and this provides an effective and simple way to construct them explicitly. On the other hand, the linear stability of such motions has been investigated in, e.g.,  \cite{Kimura, Aref10}.
 
 We shall come back to these collapses in Section \ref{subsec:filaments-collapse}, where we will use them to construct collapsing filament solutions.

\medskip

For $N\geq 4$ the $N$-vortex problem is not integrable in general, but it still exhibits some remarkable solutions.
First well-known relative equilibria  are given by the uniformly rotating configuration of $N$ identical vortices placed on the same line in a particular way, or, when $N$ is odd, 
by $N-1$ identical vortices and a central  vortex of arbitrary circulation. Other relative equilibria are  obtained by setting one vortex with same circulation $\Gamma_j=\Gamma$ at each vertex of a regular $N$ - polygon.  The motion is then a uniform rotation with constant angular velocity 
$\omega=\Gamma(N-1)/(2R^2)$, where $R$ denotes the size of the polygon.  The problem of stability of these $N$-polygons was raised by Thomson (Lord Kelvin) in the nineteenth century and has attracted much attention since then. The full answer to Thomson's question -- that  the  regular $N$-polygon is stable if and only if $N\leq 7$ -- was finally  recently provided by Kurakin and Yudovitch \cite{KY, KY2}. One can get other uniformly rotating relative equilibria by adding to a regular $N$-polygon configuration one point of arbitrary circulation $\Gamma_0$ at the center of the polygon, in which case the angular velocity is $\omega=(\Gamma(N-1)+2\Gamma_0)(2R^2)$. In particular, choosing $\Gamma_0=-\Gamma(N-1)/2$ one obtains a fixed equilibrium solution. The stability of centered $N$-polygons has been proved to hold by Cabral and Schmidt \cite{CS} for all $N$ provided the circulation $\Gamma_0$ lies in a suitable interval depending on $N$. In particular, the range of stability includes negative values of $\Gamma_0$ for $N<7$.

Finally,  we mention that other relative equilibria can be constructed in the form of  concentric rings of polygons with identical vortices. The polygons may be not regular, and the rings may have different circulations (heterogeneous rings). Relative equilibria of identical two-vortex or three-vortex rings have been determined by
Aref, Newton, Stremler,  Tokieda and  Vainchtein
  \cite{ANSTV}, and O'Neil \cite{Neil07-rings} constructed relative equilibria of heterogeneous three-vortex rings.

\medskip

The existence of collapsing solutions for $N\geq 4$ is a difficult issue. Given a set of circulations such that $\sum_{J} \Gamma_j\neq 0$ for all $J\subset \{1,\ldots,N\}$,  the set of initial configurations leading to a  collapse in finite time has zero Lebesgue measure in $\R^{2N}$ \cite{MP}.  On the other hand, in general it is not known whether 
self-similar collapses exist. Necessary  conditions for such motions are $T=0$ and $\Gamma_1\Gamma_2+\Gamma_2\Gamma_3+\Gamma_3\Gamma_1=0$, but, in contrast with the case $N=3$, there are not sufficient. Nevertheless, some specific cases have been explored:  four and five-vortex collapses were studied by Novikov and Sedov \cite{Nobikov-Sedov} and  collapses of two-vortex rings were constructed in \cite{Koiller}. Collapsing configurations for infinite point vortex lattices are constructed by O'Neil in  \cite{Neil89}.  Finally,  \cite{Neil07-rings} exhibits self-similar collapses of heterogeneous three-vortex rings.

\medskip

\section{Several filaments}\label{sec:filaments}

\subsection{The interaction of several filaments} The objective of this paragraph is to review known existence results concerning the system \eqref{syst:filaments}, which have been established by Klein, Majda and Damodaran \cite{KlMaDa} and Kenig, Ponce and Vega \cite{KePoVe}. We also mention the work by Zhakarov \cite{Zh1,Zh2}. We recall the system \eqref{syst:filaments} supplemented by initial data
$\Psi_{j,0}$:
\begin{equation*}
\begin{cases}
\displaystyle 
i\partial_t \Psi_j+\alpha_j\Gamma_j \partial_\sigma^2\Psi_j+\sum_{k\neq j} \Gamma_k \frac{\Psi_j-\Psi_k}{|\Psi_j-\Psi_k|^2}=0,\quad 1\leq j\leq N,\\
\Psi_j(0,\sigma)=\Psi_{j,0}(\sigma). 
\end{cases}
\end{equation*}

\subsubsection{Some formal properties}
\label{para:conserved}For any $N\geq 1$ and any choice of circulations $\Gamma_j\in \R$, 
the dynamics formally preserves the following quantities:

The Hamiltonian
\begin{equation*}
\frac{1}{2} \sum_{j=1}^N
\alpha_j  \Gamma_j^2 \int \left|\partial_\sigma \P_j\right|^2\,d\sigma
-\frac{1}{2}\sum_{j\neq k} \int \Gamma_j \Gamma_k  \ln \left|\P_{j}-\Psi_k\right|^2\,d\sigma;
\end{equation*}
The mean angular momentum
\begin{equation*}
\int  \sum_{j=1}^N \Gamma_j \left|\P_j\right|^2\,d\sigma,
\end{equation*}
and the mean center of vorticity
\begin{equation*}
\int  \sum_{j=1}^N \Gamma_j  \P_j\,d\sigma.
\end{equation*}
They are the three-dimensional analogs  of the notions of Hamiltonian,  angular momentum  and center of vorticity  for point vortices, which have been defined in Section \ref{Sec:point}.
There is also conservation of the quantity
\begin{equation*}
\int \sum_{j=1}^N \Gamma_j\Im m(\Psi_j \overline{\partial_\sigma \Psi_j})\,d\sigma, 
\end{equation*}
which can be interpreted as a mean momentum.

Finally, in the case where the vortex core parameters and the circulations satisfy
\begin{equation*}
\alpha_j\Gamma_j=\kappa_0\quad \forall 1\leq j\leq N
\end{equation*}
we also have conservation of the quantity
\begin{equation*}
\sum_{j\neq k}\Gamma_j \Gamma_k \int \left|\P_{j}-\Psi_k\right|^2\,d\sigma.
\end{equation*}

Nevertheless, the previous quantities may be not well-defined, not even formally. Actually, for general configurations of exactly parallel vortex filaments there are infinite. In the next section we will bypass this difficulty by the use of renormalized quantities.

\subsubsection{Some existence results}
\label{subsec:results-filaments-before}
First rigorous results on vortex filaments were obtained by Klein, Majda and Damodaran \cite{KlMaDa}  for pairs of filaments 
\begin{equation}\label{syst:pair-general}
\begin{cases}
\displaystyle 
i\partial_t \Psi_1+\alpha_1\Gamma_1 \partial_\sigma^2\Psi_1+ \Gamma_2 \frac{\Psi_1-\Psi_2}{|\Psi_1-\Psi_2|^2}=0\\
\displaystyle 
i\partial_t \Psi_2+\alpha_2\Gamma_2 \partial_\sigma^2\Psi_2- \Gamma_1 \frac{\Psi_1-\Psi_2}{|\Psi_1-\Psi_2|^2}=0.
\end{cases}
\end{equation}
It may be convenient to introduce the new variables $\Psi=\Psi_1-\Psi_2$, $\Phi=\Psi_1+\Psi_2$ and $\kappa_1=(\alpha_1\Gamma_1+\alpha_2\Gamma_2)/2$,
$\kappa_2=(\alpha_1\Gamma_1-\alpha_2\Gamma_2)/2$. Then
\begin{equation}\label{syst:pair}
\begin{cases}
\displaystyle 
i\partial_t \Phi+\kappa_1 \partial_\sigma^2\Phi+ \kappa_2 \partial_\sigma^2\Psi+(\Gamma_2-\Gamma_1)\frac{\Psi}{|\Psi|^2} =0\\
\displaystyle 
i\partial_t \Psi+\kappa_1 \partial_\sigma^2\Psi+ \kappa_2 \partial_\sigma^2\Phi+(\Gamma_1+\Gamma_2)\frac{\Psi}{|\Psi|^2} =0.
\end{cases}
\end{equation}

In \cite{KlMaDa} (see also \cite{Zh1, Zh2}) only the case of identical vortex core parameters $\alpha_1=\alpha_2$ is considered. Then
for identical circulations $\Gamma_1=\Gamma_2$ we get $\kappa_2=0$ so \eqref{syst:pair} can be decoupled into two equations for $\Phi$ and $\Psi$. The solutions have the general wave-like form
\begin{equation*}
\Psi_1(t,\sigma)=\Phi(t,\sigma)+Ae^{i(k\sigma+\omega t)},\quad \Psi_2(t,\sigma)=\Phi(t,\sigma)-Ae^{i(k\sigma+\omega t)}
\end{equation*}
where $\Phi$ is any solution to the linear Schr\"odinger equation, and where the parameters $A$, $k$ and $\omega$ satisfy a consistency equality.
On the other hand, for opposite circulations  we have $\kappa_1=0$.
By Section \ref{Sec:point}, the exactly parallel solution is given by a uniformly translating vortex pair. For nearly parallel filaments one looks for solutions with the same symmetry features, i.e.  $\Psi_2(t,\sigma)=-\overline{\Psi}_1(t,\sigma)$ for all $(t,\sigma)$. The equations \eqref{syst:pair} then yield a coupled system of equations for the real variables $\Re e(\Psi_1)=\Psi/2$ and $\Im m(\Psi_1)=-i\Phi/2$. 

From a mathematical point of view, one can also consider the case when $\alpha_1\Gamma_1=\alpha_2\Gamma_2=\kappa_0$. It follows that $\kappa_2=0$ so the equations \eqref{syst:pair} are decoupled
\begin{equation}\label{syst:pair-2}
\begin{cases}
\displaystyle 
i\partial_t \Phi+\kappa_0 \partial_\sigma^2\Phi+(\Gamma_2-\Gamma_1)\frac{\Psi}{|\Psi|^2} =0\\
\displaystyle 
i\partial_t \Psi+\kappa_0 \partial_\sigma^2\Psi+(\Gamma_1+\Gamma_2)\frac{\Psi}{|\Psi|^2} =0.
\end{cases}
\end{equation}
For identical circulations one retrieves the previous co-rotating filament pair. In the case of opposite circulations, $\Psi$ satisfies the linear Schr\"odinger equation.

\medskip

The linear stability of system \eqref{syst:pair-general}, with $\alpha_1=\alpha_2$, around the exactly parallel solution to the point vortex system 
\eqref{syst:point-vortex} is analyzed in \cite{KlMaDa}. Linear instability occurs when $\Gamma_1\Gamma_2<0$, while linear stability holds for $\Gamma_1\Gamma_2>0$. Moreover,  numerical computations predict global existence in the first case, and  self-similar collapse in finite time in the second case. So, in contrast with the $2$-dimensional setting (pairs of vortex points), finite time collapses of two nearly parallel filaments might happen. Following Zhakarov \cite{Zh1}, Klein, Majda and Damodaran suggested that such collapses
should be well described by  a class of self-similar solutions to the system \eqref{syst:pair}. Also the equations \eqref{syst:pair-2} could maybe provide a class of self-similar collapses. To our knowledge, such constructions have not been rigorously achieved yet.

\medskip

The next decisive results regarding \eqref{syst:filaments} were accomplished by Kenig, Ponce and Vega \cite{KePoVe} in 2003. The article \cite{KePoVe} is presented in the survey \cite{ponce}.  In \cite{KePoVe},  existence of vortex filaments in the form of $H^1$-perturbations of the exactly parallel solution $X_j(t)$ is investigated with the ansatz
\begin{equation}
\label{ansatz:KPV}
\Psi_j(t,\sigma)=X_j(t)+u_j(t,\sigma),\quad \text{with}\quad u_j\in C(\R,H^1(\R)).
\end{equation}
Note that, thanks to the Sobolev embedding $H^1(\R)\subset L^\infty(\R)$, the ansatz \eqref{ansatz:KPV} guarantees that the filaments are well-separated as long as the quantity  $\|u_j\|_{H^1}$ is sufficiently 
small with respect to the distance between the point vortices $X_j(t)$. As long as this remains true, the $u_j$ solve a system of Schr\"odinger equations with potentials that are Lipschitz functions in the $u_j$:
\begin{equation}
 \label{syst:pert-u-1}
\begin{cases}
\displaystyle i\partial_t u_j + \alpha_j\Gamma_j\partial_\sigma^2 u_j+\sum_{k\neq j} \Gamma_k 
\left(\frac{ X_{j}-X_k+u_{j}-u_k}{|X_{j}-X_k+u_{j}-u_k|^2}-\frac{ X_{j}-X_k}{|X_{j}-X_k|^2}\right)=0\\
\displaystyle u_j(0)=u_{j,0},\quad 1\leq j\leq N.
\end{cases}
\end{equation}
This makes it possible in \cite{KePoVe} to obtain local well-posedness for \eqref{syst:pert-u-1} for any $N$, for any set of circulations $\Gamma_j$, and for sufficiently small $\|u_j(0)\|_{H^1}$,   and at least up to times of order $\min(T^\ast,|\ln \sum_j \|u_j(0)\|_{H^1}|)$, where $T^\ast$ is the first collision time between the point vortices.

\medskip

On the other hand, Kenig, Ponce and Vega \cite{KePoVe} obtained global existence for sufficiently small $\|u_j(0)\|_{H^1}$ in the two following special cases: $N=2$ and $\Gamma_1\Gamma_2>0$, or $N=3$ and 
$(X_1(t),X_2(t),X_3(t))$ is the uniformly rotating equilateral triangle solution (see Section \ref{Sec:point}). The proof of \cite{KePoVe} is based on the use of a suitable notion of energy.
We recall that for filaments given by the ansatz \eqref{ansatz:KPV} the quantities defined in \S \ref{para:conserved} do not make sense.  This is the reason why \cite{KePoVe} introduces the renormalized quantities
\begin{equation*}\begin{split}
 \mathcal{H}&=\frac{1}{2}\sum_{j}
  \alpha_j \Gamma_j^2 \int \left|\partial_\sigma \P_j\right|^2\,d\sigma
-\frac{1}{2}\int \sum_{j\neq k} \Gamma_j \Gamma_k  \ln\left(\frac{|\P_{j}-\Psi_k|^2}{|X_{j}-X_k|^2}\right)\,d\sigma\\
 \mathcal{A}&=\int \sum_j \Gamma_j  \left(|\P_j|^2-|X_j|^2\right)\,d\sigma\\
 \mathcal{T}&=\int \sum_{j\neq k}\Gamma_j \Gamma_k  \left(|\P_{j}-\Psi_k|^2-|X_j-X_k|^2\right)\,d\sigma.\end{split}
\end{equation*}
Note that in the framework of \eqref{ansatz:KPV} we have $|\Psi_j(t,\sigma)-\Psi_k(t,\sigma)|\to |X_j(t)-X_k(t)|$ when $\sigma\to \pm \infty$.
Moreover, in view of the properties of the point vortex system, the previous quantities are still formally preserved in time provided $\alpha_j\Gamma_j=\kappa_0$  for all $1\leq j\leq N$.

Finally, we also introduce the quantity
\begin{equation*}
 \mathcal{I}(t)=\frac{1}{2}\int \sum_{j\neq k}\Gamma_j \Gamma _k  \left( \frac{|\P_{j}-\Psi_k|(t)^2}{|X_{j}-X_k|(t)^2}-1\right)\,d\sigma,
\end{equation*}
which is not necessarily constant, and we define the energy 
\begin{equation}\label{def:energy-1}\begin{split}
\mathcal E((\Psi_j(t))_j)&=\mathcal H+\mathcal I(t)\\
&=\frac{1}{2}\sum_{j}\alpha_j\Gamma_j^2
  \int \left|\partial_\sigma \P_{j}(\sigma)\right|^2\,d\sigma\\
&+\frac{1}{2}\int \sum_{j\neq k} \Gamma_j\Gamma_k  \left( -\ln\left(\frac{|\Psi_{j}(\sigma)-\Psi_{k}(\sigma)|^2}{|X_{j}-X_{k}|^2}\right)
+\frac{|\Psi_{j}-\Psi_{k}|(t)^2}{|X_{j}-X_{k}|(t)^2}-1\right)\,d\sigma
.\end{split}\end{equation}
From now on we will sometimes write
\begin{equation*}
 \mathcal{E}(t)=\mathcal E((\Psi_j(t))_j)\end{equation*}
 when not misleading.
 
We insist on the fact that, in general,  the energy is not constant in time.

Now, in the cases under consideration in \cite{KePoVe} the distances $|X_j(t)-X_k(t)|=d$ are constant in time and all the same. Therefore it turns out that 
$$\mathcal E(t)=\mathcal{H}+\frac{1}{d^2}\mathcal{T}=\mathcal E(0)$$ is constant. On the other hand, for circulations having all the same sign it was noticed in \cite{KePoVe} that the Sobolev embedding combined with the convexity inequality $(x-1)^2/2\leq x-1-\ln x\leq 10(x-1)^2$ for $x\in [3/4,5/4]$ imply coercivity for the energy:  \begin{equation}\label{ineq:coercive}
\left\| \frac{|\Psi_j-\Psi_k|^2}{|X_j-X_k|^2}-1\right\|_{\infty}\leq C\mathcal{E}((\Psi_j)_j)
\end{equation}
as long as the vortex filaments $\Psi_j$ are not too far from the straight filaments $X_j$.  Finally, it can be easily seen that 
\begin{equation}\label{ineq:coercive-2}
\mathcal{E}((X_j+u_j)_j)\leq C \sum_j \|u_j\|_{H^1}^2.
\end{equation}
Therefore in the cases of small $H^1$-perturbations considered in \cite{KePoVe} the energy is and remains small for all time, hence the vortex filaments remain well-separated for all time and global existence follows. We  stress that the vortex filaments constructed in this way remain uniformly close to the straight filaments for all time.

\subsubsection{A class of weak solutions}

To conclude this section, we mention that Lions and Majda \cite{LiMa} introduced and proved the global existence of "very weak" solutions to \eqref{syst:filaments} for any $N\geq 1$, identical circulations and identical vortex core parameters. Such solutions satisfy a weak formulation of \eqref{syst:filaments} with suitable test functions, there  are $L$-periodic with respect to $\sigma$, they belong to $ C(\R,H^1(0,L))$ and they have finite energy:
$\sup_{t\geq 0} \sum_{j\neq k} \int_0^L |\ln |\Psi_j(t)-\Psi_k(t)|\,d\sigma  <\infty$.   Actually, such a definition does not exclude the possibility of collisions, since the only information is that 
\begin{equation*}
\forall t\in \R,\quad \Psi_j(t,\sigma)\neq \Psi_k(t,\sigma)\quad \text{for a.e. } \sigma.
\end{equation*}
Moreover, uniqueness in this class is not known.

\subsection{Symmetric configurations} 
\label{subsec:filaments-symmetric}
We present now the results obtained recently in \cite{BaMi}. In all the following we will take identical vortex core parameters
$$\alpha_j=1\quad \forall 1\leq j\leq N.$$
We start with a natural remark that should in part motivate the kind of perturbations we shall consider then. As we have seen in the previous sections,  the dynamics of both systems  \eqref{syst:filaments} and  \eqref{syst:point-vortex} is complicated   when $N$ becomes large, and in fact  additional symmetry conditions are needed. In this paragraph,   we shall first describe a local in time result that extends the one in \cite{KePoVe} to energy-type spaces.  We will turn then to the case of four filaments  for which global results can be obtained without requiring too many symmetries on the perturbations. Finally, for any $N\geq 2$ we will present a global existence result for dilation-rotation-type perturbations of the regular N-polygon exactly parallel configuration. In this latter case the perturbations, at any height $\sigma$, have a regular N-polygon shape. Moreover, traveling waves will be displayed in this setting. We shall complete this subsection with another class of simple perturbations.

For a positive $\omega$ we define the energy
\begin{equation}\label{def:energy-2}
\mathcal{E}(f)= \frac{1}{2}\int |\partial_\sigma f|^2\,d\sigma+\frac{\omega}{2}\int \left(|f|^2-1-\ln |f|^2\right)\,d\sigma.
\end{equation}
The motivation for this definition will become clear in the next paragraphs.
Recurrent tools in  \cite{BaMi}, which are partially inspired by \cite{KePoVe} (in particular see \eqref{ineq:coercive} and \eqref{ineq:coercive-2}), are the following facts (see Lemmas 2.1-.2.3 in  \cite{BaMi}). 
\begin{itemize}
\item If $f$ has small energy $\mathcal E(f)$, then $f$ is close to $1$,
$$\||f|^2-1\|_{L^\infty}\leq \frac 14.$$
\item If $f$ is small in $\dot H^1$ norm and if  $t$ small enough,
$$\|e^{it\partial_\sigma^2}f-f\|_{L^\infty}\leq \frac 14.$$
\item If $\||f|^2-1\|_{L^\infty}\leq 1/4$ then we can compare the energies:
$$\mathcal E_{GP}(f)
\equiv\frac{1}{2}\|\partial_\sigma f\|_{L^2}^2+\frac{\omega}{4}\||f|^2-1\|_{L^2}^2\leq 
\mathcal E(f)\leq 5\,\mathcal E_{GP}(f).$$
\item If $f$ has small energy $\mathcal E(f)$ and $h$ is small in  $H^1$, then  the energy $\mathcal{E}(f+h)$ is finite. More precisely we have,
for absolute numerical constants $C,C'$,
\begin{equation*}\begin{split}
\mathcal{E}(f+h)&\leq C\mathcal{E}_{GP}(f+h)
\leq C'\left(1+\mathcal{E}(f)\right)\left(1+\|h\|_{H^1}^2\right),\end{split}\end{equation*}and 
\begin{equation*}
 \||f+h|-1\|_{L^\infty}<1.
\end{equation*}
\end{itemize}

\subsubsection{$N\geq 2$, a local existence result in the energy space}\label{local} As recalled in \S \ref{subsec:results-filaments-before}, local existence was proved in \cite{KePoVe} for \eqref{syst:filaments} for any $N\geq 1$ and small $H^1$-perturbations  $u_{j,0}=\Psi_{j,0}-X_{j,0}$. The relevant energy functional in that setting is $$\mathcal{E}_0=\mathcal{E}((\Psi_{j,0})_j)$$ which has been defined in \eqref{def:energy-1}. We assume next that all the circulations have the same sign, for example $\Gamma_j>0$ for all $1\leq j\leq N$. Then the energy $\mathcal{E}((\Psi_j)_j)$ satisfies the properties listed above for each function $f=|\Psi_j-\Psi_k|/|X_j-X_k|$. 
It is then natural to consider the problem of perturbations with small energy as a single assumption.  On the one hand for small initial energy we obtain that $\Psi_{0,j}$ is not too far from $X_{j,0}$ and that $ \mathcal {E}_0\leq C\sum_j \|u_{j,0}\|_{H^1}^2$.  On the other hand for $u_{j,0}$ small in  $H^1$ it is easy to see that $ \mathcal {E}_0\leq C\sum_j \|u_{j,0}\|_{H^1}^2$. Therefore the energy smallness assumption is weaker than the $H^1$ smallness assumption in general. For example, it allows for a larger geometric class of  perturbations, namely the small rotations. Indeed, take a rotation and translation type perturbation of the square configuration, 
$$\Psi_{j,0}^\eps(\sigma)=e^{i\sqrt{\eps} \varphi_0(\eps \sigma)} X_{j,0}+T^\eps(\sigma),$$
with $\varphi_0\in{H^1}$ and  $\|T^\eps\|_{ H^1}=O(\eps)$. It follows that  $\mathcal E_0=O(\eps^2)$ while $ \sum_j\|u_{j,0}\|_{H^1}^2\geq O(1).$ The idea to use the energy space  rather than $H^1$ is also reinforced by the fact that in the setting of \eqref{syst:filaments}, that we shall see later to be related in some cases to the Gross-Pitaevskii equation, the $L^2$-norm of the solutions might grow in time while the energy remains controllable.

\medskip

We next denote by $d>0$ the minimal distance between the vortices $X_j(t)$ for all time (recall that $\Gamma_j>0\,\forall j$).
In \cite{BaMi} it has been shown the following local well-posedness result:

\medskip

\emph{For an initial configuration $(u_{j,0})_j$ with small energy, there exists $T>0$ and a unique solution $u=(u_j)_j\in C([0,T],H^1(\R))^N$ to the system \eqref{syst:pert-u-1}
with 
$$\sup_{0\leq t\leq T}  \|u_j(t)\|_{H^1}\leq \|u_{j,0}\|_{H^1}+\frac{d}{4},\quad 1\leq j\leq N,$$ 
and for $T$
sufficiently small such that
$$T\big(1+\mathcal E_0+\sum_j\|u_{j,0}\|_{H^1}\big)\geq C(d,(\Gamma_j)_j).$$}
The proof is based on finding a find a fixed point in the Banach space
\begin{equation*}
B_T=\left\{ w=(w_1,\ldots,w_N)\in C\left([0,T],H^1\right)^N,\quad \sup_{0\leq t\leq T} \|w(t)\|_{H^1}
\leq \frac{d}{4}\right\} 
\end{equation*}
for the operator $A(w)=(A_j(w))_j$ defined by
\begin{equation*}\begin{split}
i\omega \int_0^t 
\sum_{k\neq j} \Gamma_k 
&\left(\frac{ X_{j}+e^{i\tau  \Gamma_j\partial_\sigma^2}u_{j,0}+w_{j}-X_k-e^{i\tau  \Gamma_k\partial_\sigma^2}u_{k,0}-w_{k}}
{|X_{j}+e^{i\tau  \Gamma_j\partial_\sigma^2}u_{j,0}+w_{j}-X_k-e^{i\tau  \Gamma_k\partial_\sigma^2}u_{k,0}-w_{k}|^2}
-\frac{ X_{j}-X_k}{|X_{j}-X_k|^2}\right)\,d\tau.
\end{split}
\end{equation*}
Then the solution is given by
$$u_{j}(t)=e^{it \Gamma_j\partial_\sigma^2}u_{j,0}+w_{j}(t).$$
Notice here that  we extend $u_j$ locally from a time $t_0$ not by a fixed point for perturbations of the initial data directly, but by a fixed point argument for small $H^1$ 
perturbations $w_j$ of the linear evolutions of the initial data. We must use crucially the fact that the deviation $e^{it\partial_\sigma^2}u_{j,0}-u_{j,0}$ can be upper-bounded in $L^\infty$ in terms of the energy at the initial time $\mathcal E_0$. This is insured from the previous study of the Gross-Pitaevskii equation (see Lemma 3 in \cite{PG}). We recall the short proof: the Fourier transform of $e^{it\partial_\sigma^2}u_{j,0}-u_{j,0}$ can 
be written as $\frac{e^{-it\xi^2}-1}{\xi}\,\xi\hat u_{j,0}(\xi)$, 
so the $L^2$ norm is bounded by $C\sqrt{t}\|\partial_\sigma u_{j,0}\|_{L^2}$ and the $\dot H^1$ norm 
is bounded by $C\|\partial_\sigma u_{j,0}\|_{L^2}$, i.e. 
$$\|e^{it\Gamma_j\partial_\sigma^2}u_{j,0}-u_{j,0}\|_{H^1}\leq C(1+\sqrt{t})
\|\partial_\sigma u_{j,0}\|_{L^2}\leq  C(1+\sqrt{t}) \mathcal{E}_0.$$

As a consequence of the local in time result we obtain that the solution $(u_j)_j$ to \eqref{syst:pert-u-1} 
exists as long as the energy $\mathcal{E}(t)$ remains small. 
Indeed note that on the one hand 
the norm $\sum_j\|u_j(t)\|_{H^1}$ can grow exponentially, but it cannot blow up as long as the energy is
sufficiently small. On the other hand, as long as the energy remains small, the filaments $\Psi_j(t,\sigma)$ remain close the the straight ones $X_j(t)$, so no collapse occurs for \eqref{syst:filaments}.

\subsubsection{$N=4$, local and global results around the square configuration} In this part we are dealing with $(X_j)_j$ the square configuration with equal circulations $\Gamma_j=1$. Let $(u_{j,0})_j\in H^1(\R)^4$ and set $\Psi_{j,0}=X_{j,0}+u_{j,0}$. We introduce the quantity
\begin{equation*}
 \tilde{\mathcal{E}_0}=\max\left\{\mathcal{E}_0;
\frac{\|u_{1,0}+u_{3,0}\|_{L^2}^2}{2}+\frac{\|u_{2,0}+u_{4,0}\|_{L^2}^2}{2}\right\}.
\end{equation*}

\medskip

In Theorem 1.2. of \cite{BaMi} it is proved that

\medskip

\emph{If $\tilde{\mathcal{E}_0}$ is small enough, there exists an absolute constant $C>0$, and there exists a time $T$,
with
$$T\geq C
\min\left\{\frac{1}{{\tilde{\mathcal{E}_0}}^{1/4}\max_{j\neq k}
\|u_{j,0}-u_{k,0}\|_{L^2}^{1/2}},\frac{1}{\tilde{{\mathcal{E}_0}}^{1/3}}\right\},$$
such that there exists 
a unique corresponding solution 
$(\Psi_j)_j$
to  \eqref{syst:filaments} on $[0,T]$, satisfying 
 $\Psi_j=X_j+u_j$, with $u_j\in C\left([0,T],H^1(\R)\right)$, and such that
\begin{equation*}
\frac 34\leq \frac{|\Psi_j(t,\sigma)-\Psi_k(t,\sigma)|}{|X_j(t)-X_k(t)|}\leq \frac 54,\quad t\in [0,T],
\quad \sigma\in \R. 
\end{equation*}
Moreover, if the initial perturbation is parallelogram-shaped, namely
\begin{equation*}
 \|u_{1,0}+u_{3,0}\|_{L^2}=\|u_{2,0}+u_{4,0}\|_{L^2} =0,
\end{equation*}
then the solution $(\Psi_j)_j$ is globally defined.}

\begin{figure}
    \includegraphics[width=7cm,height=6cm]{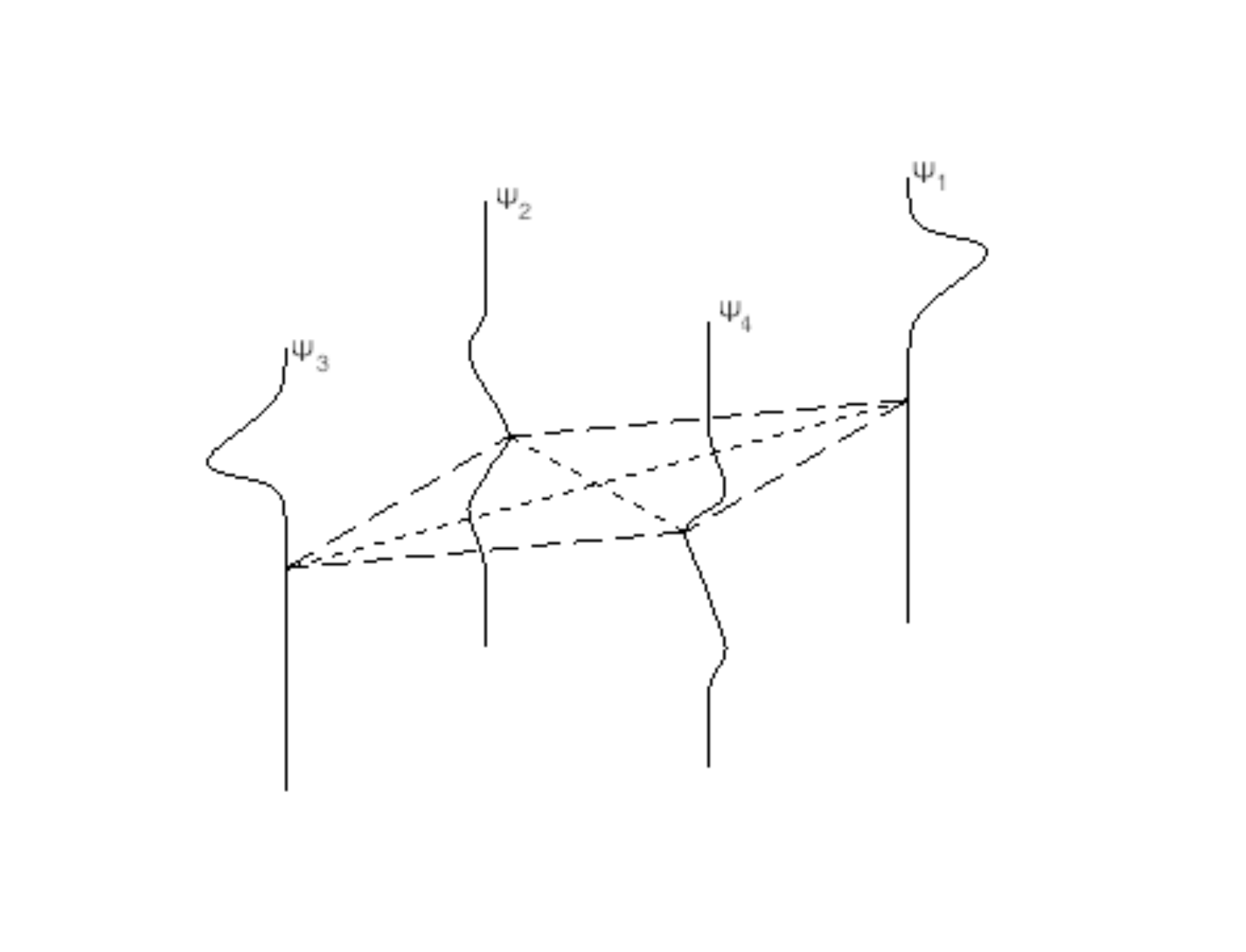} 
  \caption{$N=4$, Parallelogram-shaped perturbation around the rotating square configuration.}
    \label{Parallelogram}
\end{figure}

\bigskip

The proof goes as follows. In view of the local existence argument exposed in the previous subsection, the extension of a solution is ensured by the control of the energy. In this particular setting of the square, direct computations lead to the formula
$$ \mathcal{E}(t)=\mathcal H+\frac{1}{2}\mathcal{T}
-\mathcal{A}+
\frac{\|(u_1+u_3)(t)\|^2+ \|(u_2+u_4)(t)\|^2}{2}.$$
If $(\Psi_1,\Psi_2,\Psi_3,\Psi_4)$ is a solution of \eqref{syst:filaments} 
then $(-\Psi_3,-\Psi_4,-\Psi_1,-\Psi_2)$ is also a solution. So if the initial perturbation is a parallelogram like the one presented in Figure \ref{Parallelogram}, it will remain so at later times, hence
$$\|(u_1+u_3)(t)\|_{L^2}^2=\|(u_2+u_4)(t)\|_{L^2}^2=0$$
and global existence follows. On the other hand, for a perturbation without symmetry conditions the growth of $\|(u_1+u_3)(t)\|^2+ \|(u_2+u_4)(t)\|^2$ up to time $T$ is controlled by a lengthy computation based crucially on the particular properties of the underlying square configuration $(X_j)_j$.

Finally, let us notice that one could try the same approach for other configurations $(X_j)_j$. Explicit computations give for instance 
$$\mathcal E(t)=-\mathcal{H}+\mathcal{I}-\frac 32 \mathcal{A} + \frac 34 \left(\|u_1(t)\|_{L^2}^2
+\|(u_2+u_3)(t)\|_{L^2}^2\right),$$
for $N=3$ and $(X_j)_j$ at the
ends and the middle of the segment,
and
$$\mathcal E(t)=-\mathcal{H}+\mathcal{I}-\frac 72 \mathcal{A} + 
\frac 23 \sum_{j=1}^2\|(u_j+u_{j+2}+u_{j+4})(t)\|_{L^2}^2+\frac 34 \sum_{j=1}^3 \|(u_j+u_{j+3})(t)\|_{L^2}^2$$
for $N=6$ and $(X_j)_j$ a regular hexagon. Nevertheless, these quantities have no reason to be conserved, unless the perturbations retain the same shape as  $(X_j)_j$, which enters precisely the framework of the next paragraph.

\subsubsection{$N\geq 2$, global results for regular polygon-type configurations}  
\begin{figure}
    \includegraphics[width=7cm,height=6cm]{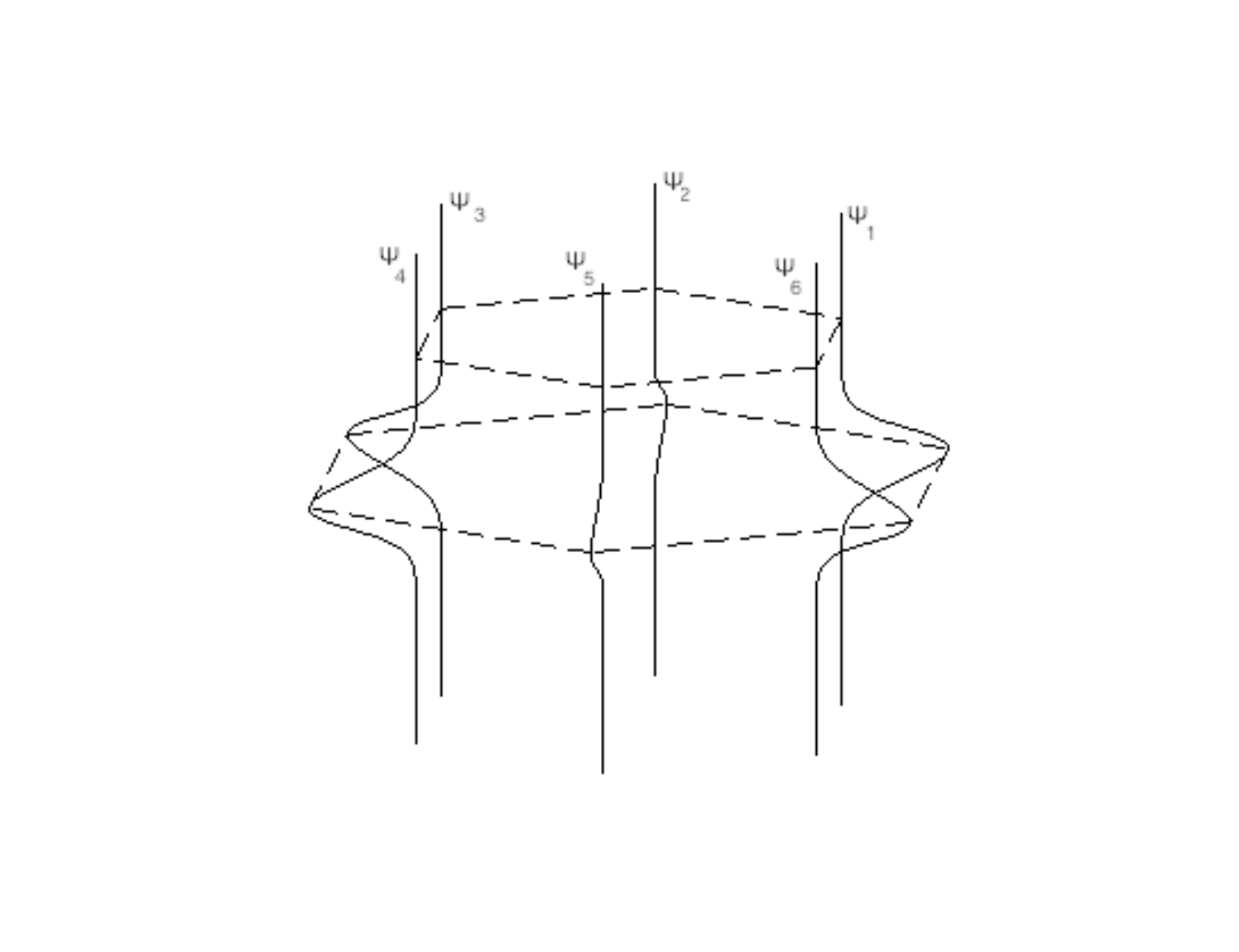} 
  \caption{$N=6$, Rotation and dilation of the polygon rotating configuration without center: $\Psi_{j}(t,\sigma)=X_{j}(t)\Phi(t,\sigma)$ for all $j$. }
    \label{Polygon}
\end{figure}

In order to get global results for a large number $N$ of filaments, in \cite{BaMi} symmetry conditions are imposed both on the straight filaments configuration and on the perturbations. More precisely, we let $(X_j)_j$ be the regular $N$-polygon configuration, with or without its center, rotating with constant velocity $\omega$ and centered at the origin. We consider  filaments $(\Psi_j)_j$ with equal circulations $\Gamma_j=1$ and of the type
\begin{equation}\label{poly}
\Psi_{j}(t,\sigma)=X_{j}(t)\Phi(t,\sigma).
\end{equation}
So in particular for any time $t$ and any height $\sigma$, the points $(\Psi_j(t,\sigma))_j$ form the same polygon as $(X_j(t))_j$ up to a rotation or/and dilation, see Figure \ref{Polygon}.

A straightforward computation shows that if $(\Psi_j)_j$ is a solution of \eqref{syst:filaments} then $\Phi$ is a solution of
\begin{equation}\label{eq:BM}
 i\partial_t \Phi+\partial_\sigma^2 \Phi+\omega \frac {\Phi}{|\Phi|^2}(1-|\Phi|^2)=0.
\end{equation}
Eq. \eqref{eq:BM} is an Hamiltonian equation, with Hamiltonian given by 
$$ \mathcal E(\Phi)= \frac{1}{2}\int |\partial_\sigma \Phi|^2\,d\sigma+\frac{\omega}{2}\int \left(|\Phi|^2-1-\ln |\Phi|^2\right)\,d\sigma.$$
In our context, as we expect that $(\Psi_j)_j$ are small nearly parallel perturbations of $(X_j)_j$, we have $|\Phi|\simeq 1$, so that  \eqref{eq:BM}  is formally similar to the already mentioned Gross-Pitaevskii 
equation
\begin{equation}\label{eq:GP}
 i\partial_t \Phi+\partial_\sigma^2 \Phi+\omega \Phi(1-|\Phi|^2)=0,
\end{equation}
with corresponding Hamiltonian  given by 
\begin{equation*}
 \mathcal E_{GP}(\Phi)= \frac{1}{2}\int |\partial_\sigma \Phi|^2+\frac{\omega}{4}
\int \left(|\Phi|^2-1\right)^2. 
\end{equation*}
This formal comparison turns out to be valid: as we have reminded at the beginning of this subsection, when $|\Phi|\simeq 1$ the energies $\mathcal{E}(\Phi)$ and $\mathcal{E}_{GP}(\Phi)$ are comparable. So our objective is now to solve globally the Cauchy problem for the equation \eqref{eq:BM} in the energy space, in the same spirit as previous works on \eqref{eq:GP}, see \cite{PG}.

\medskip

We next assume  that $\omega>0$, which is the case as soon as the circulation $\Gamma_0$ of the central vortex satisfies $N-1+2\Gamma_0>0$ (see Section \ref{Sec:point}).
Theorem 1.1 of \cite{BaMi} asserts that

\medskip

\emph{For an initial data
$$\Psi_{j,0}(\sigma)=X_{j}(0)\Phi_0(\sigma),$$
with $(X_j)_j$ a regular polygon configuration, with or without its center, rotating with constant positive speed $\omega$, and $\mathcal E(\Phi_0)$ small, \eqref{syst:filaments} has a unique global solution of type \eqref{poly} such that $\mathcal E(\Phi(t))=\mathcal E(\Phi_0)$ and $\Phi-\Phi_0\in C(\R,H^1(\R))$. 
In particular, if $\Phi_0(\sigma)\overset{|\sigma|\rightarrow\infty}{\longrightarrow}1$ then for all time $\Psi_j(t,\sigma)\overset{|\sigma|\rightarrow\infty}{\longrightarrow} X_j(t)$, so the filaments $(\Psi_j)_j$ remain nearly parallel.}

\medskip

The proof  is based on the transposition of the arguments presented in \S \ref{local} to the setting of \eqref{eq:BM}.  We notice that 
$$\mathcal E=N\mathcal E(\Phi).$$
Since  $\mathcal E(\Phi)$ is conserved it follows that the corresponding local existence result for \eqref{eq:BM} can be iterated and global existence is achieved. 

In the third part of \cite{BaMi} subsonic traveling waves for \eqref{syst:filaments} are constructed, still under the condition $\omega>0$. More precisely, for $c<\sqrt{2\omega}$ sufficiently close to $\sqrt{2\omega}$ there exist solutions of \eqref{syst:filaments} of the shape
$$\Psi_j(t,\sigma)=X_j(t)v(\sigma+ct),$$
with $v(\sigma)$ a smooth function of small energy, of even increasing modulus, exponentially increasing to $1$ as $\sigma$ tends to $\infty$. Moreover, 
$$v(\sigma)\overset{\sigma\rightarrow\pm\infty}{\longrightarrow} e^{i\theta_\pm},\quad|\theta_+-\theta_-|\leq C\sqrt{2\omega-c^2}.$$

These traveling waves are the equivalents objects of the grey solitons for the 1-D Gross-Pitaevskii equation. The proof of the existence follows the approach in \cite{Ma} and \cite{Gr}. The equation for the profile $\eta=1-|v|^2$ 
$$\eta'=\Big(-(c^2-4\omega)\eta^2+4\omega\big( (\eta-1)\ln(1-\eta)-\eta\big) \Big)^{1/2}
,$$
 is a little more complicated than the one for the Gross-Pitaevskii equation which can be solved explicitly. It is  integrated in \cite{BaMi} using the smallness assumption on $c-\sqrt{2\omega}$. 

As long as the energy is small, the modulus of $\Phi(t,\sigma)$ stays close to $1$ and the nonlinearity in \eqref{eq:BM} enters the frameworks of \cite{ZL} (see also \cite{BeGrSa}), so these traveling waves are orbitally stable. Also, the results in \cite{Ma} are valid and insure that there are no non-trivial supersonic traveling waves. 

Finally, by reusing an argument from \cite{KePoVe}, we can transform these traveling waves in the following way. System  \eqref{syst:filaments} is Galilean invariant, which means that for any real $\nu$ and any solution $(\Psi_j)_j$, the family $\Psi_{j,\nu}(t,\sigma)=e^{-it\nu^2+i\nu \sigma}\,\Psi_j(t,\sigma-2t\nu)$ is
  also a solution. By choosing $\nu=\sqrt{\omega}$ we get as another solution of \eqref{syst:filaments}
 $$\tilde\Psi_{j}(t,\sigma)=e^{i\sqrt{\omega} \sigma+ i\frac{2\pi j}{N}}\,v(\sigma+t(c-2\sqrt{\omega})),$$
which represents a stationary $(\theta_+-\theta_-)$-twisted $N$-helix filament configuration with some localized perturbation traveling in time on each of its filaments.

\subsubsection{$N\geq 2$, shifted perturbations} \label{shift} As noticed in \cite{BaMi}, if the perturbations are considered of the shifted form 
$$\Psi_j(t,\sigma)=X_j(t)+u(t,\sigma),$$ 
for any 
$(X_j)_j$ with $\Gamma_j$ all the same, it follows that $u$ must solve the linear Schr\"odiger equation.  For an initial data in $L^1$, by using the dispersion inequality for the linear Schr\"odinger equation it follows that the filaments remains separate globaly in time and the perturbations spread along the straight configuration $(X_j)_j$. By considering less regular initial data one can get examples of perturbations $u$ decaying at infinity that lead to a $L^\infty$ dispersive blow-up for the linear Schr\"odinger equation. For instance the homogeneous data $|x|^{-p}$ with 
$0<p<1$ yields a self-similar linear Schr\"odinger solution, smooth for positive times (\cite{CaWe}). So this solution leads to solutions of \eqref{syst:filaments} blowing-up in $L^\infty$ in finite time at height $\sigma=0$. Recently an example of initial data in $L^2$ but not in $L^1$ that provides dispersive blow-up was given in \cite{BoSa}. More precisely, the solution with initial data $e^{i|x|^2}/(1+|x|^2)^m$ with $1/2<m\leq 1$ blows-up in finite time at one point. This makes a second example of solutions of \eqref{syst:filaments} that blow up in finite time at a certain height but for which no collision occur.

\subsection{Collapses of filaments}
\label{subsec:filaments-collapse}
In the last part of this paper we describe configurations of nearly parallel filaments evolving towards collision in finite time. Most of them are new results. They are based on perturbations of type \eqref{poly} of exactly parallel filament configurations $(X_j)_j$. Therefore the collapse for the filaments $(\Psi_j)_j$ solutions of \eqref{syst:filaments} is linked to solutions of \eqref{eq:BM}, with modulus initially close to one, that vanish at least at one point in finite time. So a pointwise control of solutions of \eqref{eq:BM} is needed, which is quite unusual in the study of the Schr\"odinger equation. 

Let us first notice that if $(\Psi_{j,0})_j$ leads to a collision in finite 
time, the shifted perturbations introduced in \S\ref{shift},
$$\tilde\Psi_{j,0}(\sigma)=\Psi_{j,0}+u_0(\sigma),$$
with $u_0$ in $H^1$ for instance, yield a collapse of the same kind. So 
collapses are not isolated phenomena.

\subsubsection{Gaussian collapses about a stationary regular polygon configuration}
\begin{figure}
    \includegraphics[width=7cm,height=6cm]{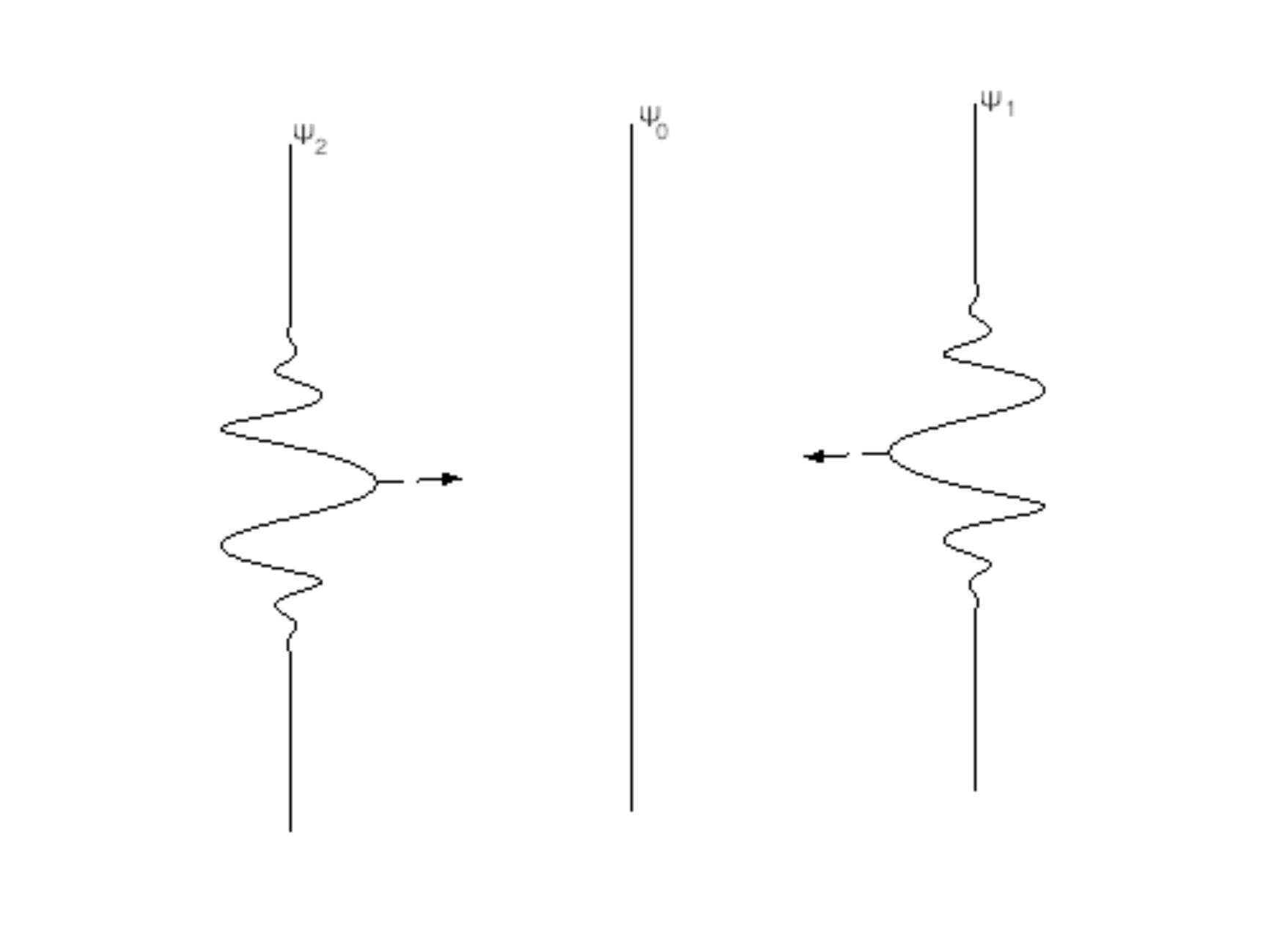} 
  \caption{$N=2$, Gaussian perturbation around the fixed equilibrium solution $(X_1,X_0=0,X_2=-X_1)$, leading to collision in finite time. }
    \label{Gaussian-collapse}
\end{figure}

We start by recalling the single collapse situation described in \cite{BaMi}. We let $(X_j)_j$ be the stationary regular $N$-polygon configuration, with circulation $\Gamma_j=1$ in the $N$ vertices of the polygon and circulation $-(N-1)/2$ in its center (see \S \ref{section:vortexpoints}). Note that this implies $N\geq 2$. Since the angular velocity is $\omega=0$, it follows that $\Phi$ is a solution of the linear Schr\"odinger equation. The evolution of a real Gaussian $e^{-\sigma^2}$ by the linear Schr\"odinger equation can be computed explicitly  
$$\frac{e^{-\frac{\sigma^2}{1+4it}}}{\sqrt{1+4it}}.$$
Therefore we are insured that the non-vanishing initial condition 
$$\Psi_{j,0}(\sigma)=X_j(0)\left(1-\frac{e^{-\frac{\sigma^2}{1-4i}}}{\sqrt{1-4i}}\right)$$
yields a solution $(\Psi_j)_j$ for the system  \eqref{syst:filaments},
$$\Psi_j(t,\sigma)=X_j(t)\left(1-\frac{e^{-\frac{\sigma^2}{1-4i(1-t)}}}{\sqrt{1-4i(1-t)}}\right),$$
 with $\Psi_j-X_j\in C\left(\R,H^1(\R)\right)$, 
such that all filament collide together at time $t=1$ at $\sigma =0$. 
Figure \ref{Gaussian-collapse} displays such initial perturbations for $N=2$ around the stationary configuration $(X_1,X_2=-X_1)$  with center at $X_0=0$.

\medskip
As a new observation, we can exhibit a similar scenario of collision for special perturbations of the previous data, namely for all initial data of type
$$\Psi_{j,0}(\sigma)=X_j(0)\left(1-\frac{e^{-\frac{\sigma^2}{1-4i}}}{\sqrt{1-4i}}+e^{i\partial_\sigma^2}u_0(\sigma)\right),$$
with $u_0$ small enough in $L^1$ and such that 
$$\int e^{i\frac{y^2}{8}}u_0(y)dy=0.$$
This is due to the fact that  
$$\Psi_j(t,\sigma)=X_j(t)\left(1-\frac{e^{-\frac{\sigma^2}{1-4i(1-t)}}}{\sqrt{1-4i(1-t)}}+e^{i(t+1)\partial_\sigma^2}u_0(\sigma)\right)$$
is then a solution of the system  \eqref{syst:filaments}. In view of the dispersion inequality and the explicit formula for a linear Schr\"odinger solution, we have
$$\left\|e^{i(t+1)\partial_\sigma^2}u_0\right\|_{L^\infty}\leq \frac{C}{\sqrt{t+1}}\|u_0\|_{L^1},\quad e^{i2\partial_\sigma^2}u_0(0)=C\int e^{i\frac{y^2}{8}}u_0(y)dy.$$
Therefore we are insured that at initial time the filaments $(\Psi_j)_j$ are separated and they collide at time $t=1$ at $\sigma=0$.

\subsubsection{Self-similar collapses about a stationary regular polygon configuration}Again we consider $(X_j)_j$ to be the stationary regular polygon configuration, with circulations $\Gamma_j=1$ in the vertices of the polygon and circulation $-(N-1)/2$ in its center. This time we take advantage of the explicit self-similar solution of the linear Schr\"odinger equation exhibited by Cazenave and Weissler \cite{CaWe}. The initial data considered there is $\psi(x)=|x|^{-p}$ with 
$0<p<1$. 
By combining (3.2), Proposition 3.3. and Corollary 3.4. b) of  \cite{CaWe} it follows that 
$$(4it)^\frac p2e^{it\partial_\sigma^2}\psi (0)\in\mathbb R,\quad \left\|(4it)^\frac p2e^{it\partial_\sigma^2}\psi\right\|_{L^\infty}=(4it)^\frac p2e^{it\partial_\sigma^2}\psi (0)=\frac{\Gamma\left(\frac{1-p}{2}\right)}{\Gamma\left(\frac 12\right)}.$$  
In particular we obtain the existence of two times $0<t_2<t_1$ such that
$$\left\|i^\frac p2e^{it_1\partial_\sigma^2}\psi\right\|_{L^\infty}<\frac 12,\quad i^\frac p2e^{it_2\partial_\sigma^2}\psi (0)=1.$$
Therefore we consider as initial data the non-vanishing function 
$$\Psi_{j,0}(\sigma)=X_j(0)\left(1-\overline{i^\frac p2e^{it_1\partial_\sigma^2}\psi (\sigma)}\right).$$
It yields as a solution for \eqref{syst:filaments}
$$\Psi_j(t,\sigma)=X_j(0)\left(1-\overline{i^\frac p2e^{i(t_1-t)\partial_\sigma^2}\psi (\sigma)}\right),$$
which vanishes at time $t=t_1-t_2$ at $\sigma=0$. Moreover, Proposition 3.7. in \cite{CaWe} insures that $e^{it\partial_\sigma^2}\psi$ belongs to some $L^r_\sigma$ space, so we have indeed $\Psi_j(t,\sigma)\overset{|\sigma|\rightarrow\infty}{\longrightarrow} X_j(t)$.

\subsubsection{Self-similar collapse around a three-vortex collapse}

The purpose of this paragraph is to investigate the behavior of vortex filaments around a self-similar  collapsing solution to the three-vortex problem.
We consider the one constructed in, e.g.,  \cite{Aref79} or \cite{MP}. The initial configuration is
\begin{equation*}
 \Gamma_1=\Gamma_2=2,\quad \Gamma_3=-1,\quad X_1(0)=-1,\quad X_2(0)=1,\quad X_3(0)=1+i\sqrt{2}.
\end{equation*}
 The sufficient collapse conditions $\Gamma_1\Gamma_2+\Gamma_2\Gamma_3+\Gamma_3\Gamma_1=0$ and $T=0$ are both satisfied.
  Then if we denote by
$
 c=-\frac{1}{3}(1+i\sqrt{2})
$
the corresponding center of vorticity, the configuration satisfies
\begin{equation*}
\dot{X}_j(t)=i\left(X_j(t)-c\right)\omega(t),\quad 1\leq j\leq 3,
\end{equation*}
where 
\begin{equation*}
 \omega(t)=\frac{\overline{\omega}}{1-\frac{t}{\tau}},\quad \overline{\omega}=
a+ib= \frac{5}{6}+i\frac{\sqrt{2}}{6},\quad \tau=\frac{1}{2b}=\frac{3}{\sqrt{2}}.
\end{equation*}
The solution  is  explicitly given by
\begin{equation}\label{sol:three}
 X_j(t)-c=\left(X_j(0)-c\right)\sqrt{1-\frac{t}{\tau}}\exp\left(-ia \tau
 \ln\left(1-\frac{t}{\tau}\right)\right),\quad t\in[0,\tau).
\end{equation}

\medskip

\begin{figure}
    \includegraphics[width=7cm,height=5cm]{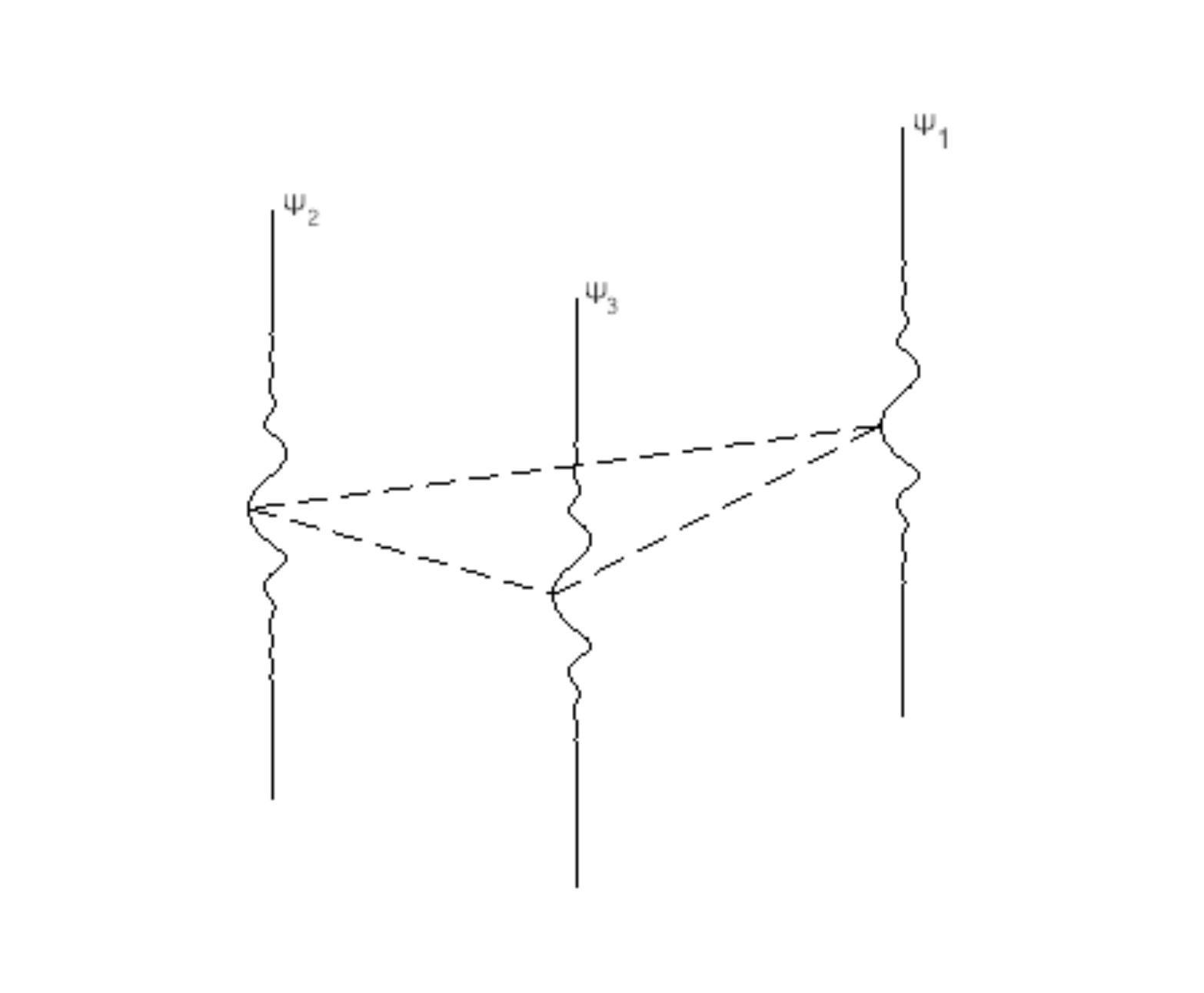} 
  \caption{$N=3$, Shifted perturbation around a self-similar collapsing triangle configuration:  $\Psi_{j,0}(\sigma)=X_j(0)+u_0(\sigma)$ for all $j$.}
    \label{Triangle-collapse}
\end{figure}

As already mentioned, since the exact parallel solution $(X_1(t),X_2(t),X_3(t))$ collapses in finite time then any shifted perturbations of the 
form 
$\Psi_{j,0}(\sigma)=X_j(0)+u_0(\sigma)$, for all $1\leq j\leq 3$, with $u_0\in H^1$, leads to collision at time $\tau$ (see Figure \ref{Triangle-collapse}).

\medskip

We next look for other kinds of less trivial perturbations around the triangle collapse. As we did before in Section \ref{subsec:filaments-symmetric}, we look for a filament solution under the form
\begin{equation*}
 \Psi_j(t,\sigma)-c=\left(X_j(t)-c\right)\Phi(t,\sigma),\quad \text{with} \: |\Phi(0,\sigma)|\simeq 1.
\end{equation*}

From now on we will assume that 
\begin{equation*}
 \alpha_j\Gamma_j=\kappa_0, \quad1\leq  j\leq 3. 
\end{equation*}
Then the equation for $\Phi$ is
\begin{equation*}
 i\partial_t \Phi+\kappa_0 \partial_\sigma^2 \Phi + 
 \frac{-i\dot{X}_j}{X_j-c}\frac{\Phi}{|\Phi|^2}\left(1-|\Phi|^2\right)=0,
\end{equation*}
so finally
\begin{equation}\label{eq:collapse-triangle}
i\partial_t \Phi+\kappa_0 \partial_\sigma^2 \Phi + 
 \omega(t)\frac{\Phi}{|\Phi|^2}\left(1-|\Phi|^2\right)=0.
\end{equation}
In the case of an equilateral uniformly rotating triangle we would have $\omega(t)\equiv \omega\in \R$ real and constant and we would retrieve Eq. \eqref{eq:BM}.

We next investigate the existence of self-similar solutions to \eqref{eq:collapse-triangle}. We seek for 
a solution in the form
\begin{equation}\label{ansatz:phi}
 \Phi(t,\sigma)=r(t)\exp(i\beta(t))\exp\left(i \frac{\sigma^2}{\gamma(t)}\right),\quad r(t)\geq 0,\quad \gamma(t)\in \R,
 \end{equation}
 where $r(0)=r_0\simeq 1$ satisfies $r_0<1$.
 
 \medskip

Plugging \eqref{ansatz:phi} into \eqref{eq:collapse-triangle} we find the system  (recall that $\tau b=1/2$)

\begin{equation}
 \label{syst:self-sim}
\begin{cases}
\dsp \dot{r}+\frac{2 \kappa_0 r}{\gamma} +\frac{1}{r}(1-r^2)\frac{1}{2(\tau-t)}=0\\
 \dsp \dot{\gamma}-4\kappa_0 =0\\
\dsp -r\dot{\beta}+\frac{1}{r}(1-r^2)\frac{a\tau}{\tau-t}=0.
\end{cases}
\end{equation}
By setting $$\gamma(t)=4\kappa_0(t-\alpha\tau),$$
where $\alpha>0$ is a parameter, we can find an explicit formula for $r$:
\begin{equation}\label{eq:r}
 r^2(t)=\frac{r_0^2-\frac{\alpha}{2}}{(1-\frac{t}{\tau})(1-\frac{t}{\alpha\tau})}+\frac{\alpha}{2}\frac{1-\frac{t}{\alpha \tau}}{1-\frac{t}{\tau}},\quad t\in  I,
\end{equation}
where $I$ is the largest interval such that $I\subset [0,\min(\tau,\alpha\tau))$ and $r>0$ on $I$.

\medskip

We rewrite \eqref{eq:r} as
\begin{equation}\label{eq:r-2}r^2(t)=\frac{1}{2\alpha \tau^2}\frac{1}{(1-\frac{t}{\tau})(1-\frac{t}{\alpha\tau})}(t^2-2\alpha \tau t+2\alpha \tau^2 r_0^2)\end{equation} and observe that the discriminant in the numerator is negative for $\alpha<2r_0^2$ and positive for $\alpha>2r_0^2$. We are led to the study of different cases which exhibit different kinds of motion. From now on $C$ will denote a positive constant depending only on $\tau$ and $\alpha$, which can possibly change
from one line to another.

\medskip

\noindent \textbullet \, We have $0<\alpha < 2r_0^2$. It follows that $ I=  [0,\min(\tau,\alpha \tau))$ and
$$ \frac{1}{C\sqrt{\left(1-\frac{t}{\tau}\right)\left(1-\frac{t}{\alpha\tau}\right)}}\leq  r(t)\leq  \frac{C}{\sqrt{\left(1-\frac{t}{\tau}\right)\left(1-\frac{t}{\alpha\tau}\right)}}\,\quad t\in  [0,\min(\tau,\alpha \tau)).$$
After finding $\beta$ in \eqref{syst:self-sim}, and taking into account \eqref{sol:three} we obtain that the filament solution constructed this way satisfies:
\begin{equation*}\begin{split}
 &\frac{1}{C\sqrt{1-\frac{t}{\alpha\tau}}}\leq  |(\Psi_j-c)(t,\sigma)|\leq  \frac{C}{\sqrt{1-\frac{t}{\alpha\tau}}},\\
&  \frac{1}{C\sqrt{1-\frac{t}{\alpha\tau}}}\leq  |(\Psi_j-\Psi_k)(t,\sigma)|\leq \frac{C}{\sqrt{1-\frac{t}{\alpha\tau}}}, \quad \sigma\in \R, \quad t\in [0,\min(\tau,\alpha\tau)). \end{split}\end{equation*}
If $\alpha\leq 1$, this means that a blow-up of the filament solution at each level $\sigma$, without collapse, takes place at time $t=\alpha\tau$ (strictly before the occurrence of the three-vortex collapse if $\alpha<1$ and exactly at the same time if $\alpha=1$). If $\alpha>1$, no blow-up nor collapse occurs for the filaments at the  time $t=\tau$ of the three-vortex collapse. Actually the modulus of the filaments remain regular up to $t=\tau$ and only the angular velocity becomes singular at time $\tau$.

\medskip

\noindent \textbullet \, We have $2r_0^2\leq\alpha<\frac{1}{2(1-r_0^2)}$. We formulate \eqref{eq:r-2} as
$$r^2(t)=\frac{1}{2\alpha \tau^2}\frac{1}{(1-\frac{t}{\tau})(1-\frac{t}{\alpha\tau})}(t-t_\ast)(t-t^\ast), \quad t\in I,$$
with $t_\ast=\alpha \tau \left[1+\sqrt{1-\frac{2r_0^2}{\alpha}}\right]$ and $t^\ast=\alpha \tau \left[1-\sqrt{1-\frac{2r_0^2}{\alpha}}\right]$. Since both $t_\ast>\tau$ and $t^\ast>\tau$ we get
the same behavior as for the case $1<\alpha<2r_0^2$ for the filaments.

\medskip

\noindent \textbullet \, We have $\frac{1}{2(1-r_0^2)}<\alpha$. Then $t_\ast>\tau$ but $0<t^\ast<\tau$; therefore, $I=[0,t^\ast)$ is strictly contained in $[0,\tau)$. 
We conclude that the filament solution evolves toward the following self-similar collision in finite time: it exists up to time $t^\ast<\tau$, collapses at each level $\sigma$ at time $t=t^\ast$ and satisfies 
\begin{equation*}
|\Psi_j(t,\sigma)-c|\sim C\sqrt{1-\frac{t}{t^\ast}�},\quad  | (\Psi_j-\Psi_k)(t,\sigma)|\sim C\sqrt{1-\frac{t}{t^\ast}},\quad t\to t^\ast,\quad \sigma\in \R.
\end{equation*}
In particular the collision between the filaments takes place before the collision between the point vortices.

\medskip

\begin{remark} When $r_0>1$ the last argument does not apply and we are not able to construct a collapse in that way.
\end{remark}

\begin{remark} In the case of a uniformly rotating polygon (see \S \ref{subsec:filaments-symmetric}) we have $\omega(t)\equiv \omega$ is constant 
in time, real, so that the analog of \eqref{syst:self-sim} is
\begin{equation*}
 \begin{cases}
  \dsp \dot{r}+\frac{r}{2t}=0\\
\dsp \dot{\beta}=\frac{\omega}{r^2}(1-r^2)
 \end{cases}
\end{equation*}
 Thus the corresponding filament solution is given by
\begin{equation*}
 \Psi_j(t,\sigma)=X_j(t)\frac{u_0}{\sqrt{1-\frac{t}{\tau}}}\exp\left( i \,\left(\frac{\omega \tau^2}{2}(1-\frac{t}{\tau})^2+\omega\tau(1-\frac{t}{\tau})\right)\right) \exp\left(\frac{i\sigma^2}{4\kappa_0 (t-\tau)}\right) 
\end{equation*}
for some $u_0\in \mathbb{S}^1$ and $\tau\in \R$. It blows up at time $t=\tau$ but does not collide.
\end{remark}

\end{document}